\documentclass[12pt]{article}
\usepackage{amsmath}
\usepackage{fullpage}

\usepackage[sort&compress,numbers]{natbib}
\usepackage{bbm}

\usepackage[utf8]{inputenc}
\usepackage[english]{babel}
\usepackage{amsmath}

\usepackage{amsfonts}
\usepackage{listings,xcolor,caption, mathtools, wrapfig}
\usepackage{subfig}
\usepackage{graphicx}
\usepackage{multicol}
\usepackage{amssymb,graphicx,enumerate,makecell}
\usepackage{pdfpages}
\usepackage{hyperref}
\usepackage{eurosym}

\usepackage{tcolorbox}
\usepackage{multirow}
\usepackage[ruled,vlined]{algorithm2e}

\usepackage[normalem]{ulem} 
\usepackage{siunitx} 
\usepackage{booktabs}
\usepackage[mathscr]{euscript}

\newtheorem{remark}{Remark}

\newcommand{\nbOut}{K}

\usepackage{color,ulem}
\definecolor{myorange}{rgb}{0.9568,0.4941,0.1961}
\definecolor{myred}{rgb}{0.9098,0.1294,0.2078}
\definecolor{myblue}{rgb}{0.0352,0.4981,0.6509}
\definecolor{myhyperblue}{rgb}{0.1607,0.3922,0.9}
\definecolor{mygreen}{rgb}{0.2235,0.6353,0.2588}
\definecolor{mygrey}{rgb}{0.3,0.3,0.3}

\newcommand{\kspace}[1]{}

\newcommand{\curves}[1]{}

\newcommand{\new}[1]{{\color{cyan}#1}}

\title{Duct boundary conditions for incompressible fluid flows:  finite element
    discretizations and
    parameter estimation in coronary blood flow
}

\author{Jerem\'ias Garay$^{1,2}$, David Nolte$^1$,
    Crist\'obal Bertoglio$^1$\footnote{Corresponding author: \texttt{c.a.bertoglio@rug.nl}}\\
    $^1$\small{Bernoulli Institute, University of Groningen, The Netherlands} \\
    $^2$\small{Department of Mechanical Engineering, Pontificia Universidad Cat\'olica, Chile}
}

\renewcommand{\vec}{\textbf}

\begin{document}

\maketitle


\begin{abstract}
    3D-0D coupled flow models are widely used across many application fields but remain
    challenging to solve. Implicit coupling introduces non-local terms, whereas explicit
    coupling results in only conditionally stable schemes. Furthermore, incorporating
    inertial effects alongside viscous resistance enlarges the parameter space, making
    calibration more difficult.

    In this work, we propose a new type of boundary condition based on the method of
    asymptotic partial decomposition of a domain (MAPDD), which we denote as the
    \textit{Duct Boundary Condition} (DuBC). This approach enables the
    incorporation of geometrically reduced domains as a boundary term with only local
    coupling in the implicit case. Moreover, the DuBC accounts for both viscous and
    inertial effects simultaneously using a single physical parameter. Additionally, we
    derive a fractional-step time-marching scheme including the DuBC. We demonstrate the
    features of the DuBC in coronary artery blood flow simulations, including sequential
    parameter estimation from noisy velocity data.


\end{abstract}



\section{Introduction }
\label{sec:introduction}

Vascular blood flow simulations of large anatomical portions are computationally
prohibitive since complex, and large geometries are usually involved. Therefore, the
usual approach is to solve the Navier-Stokes equations in region of interest,
introducing a reduced-order model obtained from geometrical assumptions as boundary
conditions to represent the remainder of the vasculature.

The typical strategy is to use \textit{0D models}, where a vessel network can be
expressed in terms of resistances, compliances, and inertance of different anatomical
subregions \cite{windkesselartery}. All these models have in common that the reduction
step is made before the coupling with the 3D geometry. This results in a non-local
coupling among the degrees of freedom on the boundary, which leads to intricate linear
algebra problems if discretized implicitly in time. However, when coupled explicitly,
instabilities may arise \cite{bertoglio2013fractional, grandmont2021existence}. For
those reasons, 0D models are challenging to use in an inverse problem setting where
the number of forward problem solutions is large and robustness to a wide range of
parameter values is required.

The method of asymptotic partial decomposition of a domain (MAPDD)
\cite{panasenko1998method,blanc1999asymptotic} is a strategy that can be used for
geometrical multiscale flow simulations, i.e., when the domain of interest contains
different levels of characteristic sizes, such as the vascular networks mentioned above.

In this approach, the vascular network is modeled as a combination of two types of
domain regions. The first are the larger regions, called \emph{junctions}, where the
blood flow is fully resolved. The second are smaller regions, represented as thin
cylindrical structures called \emph{ducts}, where the flow description is simplified.
Because the reduction is performed by selecting an appropriate subspace solution
within the region where the geometric assumption holds, the resulting formulation is,
by construction, well-posed. 

The MAPDD was first presented in \cite{panasenko1998method} using steady-state Stokes
equations in the junctions and assuming the flow inside the ducts modeled as a
Poiseuille flow, i.e., a parabolic profile for the velocity with axial symmetry and
driven by a constant pressure gradient. Later on, a generalization was made allowing
for time-dependent flows in \cite{bertoglio2019junction,bertoglio2021reconstruction}.

In this work, we formulate MAPDD directly as a boundary condition, since it only
requires the definition of a single parameter: the length of the ``extrusion'' where a
transient Stokes problem with an arbitrary time-varying pressure gradient is defined,
thanks to the constraint that the fluid flows parallel to the straight ducts. This
``virtual length'' allows us to parametrize both viscous and inertial effects
simultaneously, without assuming any specific velocity profile shape as in classical
0D models, and is therefore general for any outlet shape. We denote this MAPDD-based
boundary condition as the \textit{Duct Boundary Condition} (DuBC).

A main contribution of our work is the development and analysis of a fractional-step formulation,
whereas previous works solved MAPDD using a monolithic velocity-pressure coupling.
This approach enables a significant reduction in computational cost and allows for
efficient parameter estimation, which is an important need, e.g., in patient-specific
assessment of hemodynamic conditions \cite{nolte2022inverse}.

The rest of this article is structured as follows. In Section \ref{sec:forward} the
mathematical method is introduced followed by the description and solution of the
forward problem using 
the DuBC. In Section \ref{sec:CT_algo}, the fractional step
formulation is presented and analyzed. In Section \ref{sec:inv}, we describe and
present results of estimating
the extensions' lengths from MRI-like velocity measurements. Finally, conclusions are
given in Section \ref{sec:conclusion}.

\section{Geometric multiscale fluid flow modeling} \label{sec:forward}

\subsection{Full domain model}
\label{sec:reference_solution}

\begin{figure}[!hbtp]
    \centering
    \subfloat[ ]{
        \includegraphics[trim=200 75 100 0, clip,
    width=0.45\textwidth]{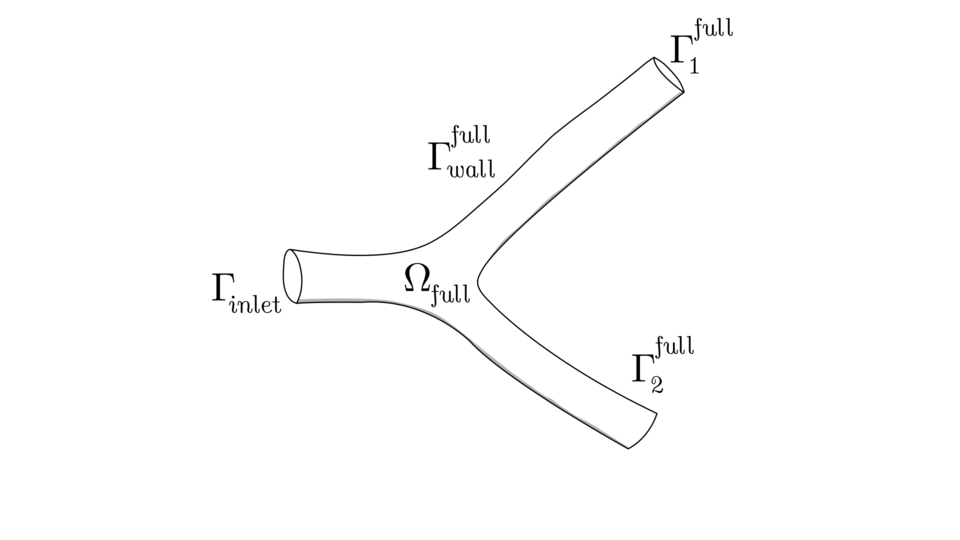} \label{fig:duct_mesh-full}}
    \subfloat[ ]{
        \includegraphics[trim=200 75 100 0, clip,
    width=0.45\textwidth]{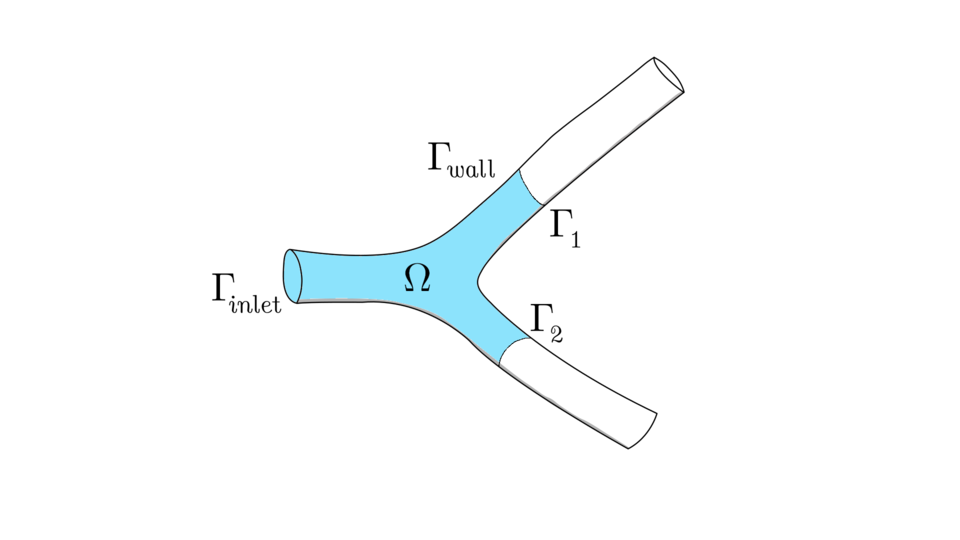} \label{fig:duct_mesh-cut}}
    \caption{Sample domains with $K=2$. (a) Full geometry; (b) Reduced geometry.}
    \label{fig:duct_mesh}
\end{figure}

Consider the incompressible Navier-Stokes equations, in a ramified domain
$\Omega_{full}$ as depicted in Figure \ref{fig:duct_mesh-full}, where its boundary is
split into an inflow $\Gamma_{inlet}$, a wall $\Gamma^{full}_{wall}$, and $K$ outlets
$\Gamma^{full}_{1},\dots,\Gamma^{full}_{K}$.

The full order model then  reads: Find $\vec{u}:\Omega \times [0,T] \to \mathbb{R}^3$,
$p:\Omega \times [0,T] \to \mathbb{R}$ such that
\begin{subequations}
    \label{eq:NS_full}
    \begin{align}
        \rho \frac{\partial \vec{u}}{\partial t} + \rho \left ( \vec{u}  \cdot \nabla
        \right ) \vec{u} - \mu \Delta \vec{u} + \nabla p  =  0 \quad \text{in} \quad
        \Omega_{full} \times (0,T] \\
        \nabla \cdot \vec{u}  =  0 \quad \text{in} \quad \Omega_{full} \times (0,T] \\
        \vec{u}  =  \vec{u}_{inlet}(\vec{x},t) \quad \text{on} \quad \Gamma_{inlet}
        \times (0,T] \\
        \mu \frac{\partial \vec{u}}{\partial \vec{n}} - p\vec{n} = \vec{0} \quad \text{on}
        \quad \Gamma^{full}_{K} \times (0,T], \quad k=1,...,K \\
        \vec{u} = \vec 0 \quad \text{on} \quad \Gamma^{full}_{wall} \times (0,T] \\
        \vec{u}(\vec{x},0) = \vec 0 \quad \text{in} \quad \Omega_{full}
    \end{align}
\end{subequations}
where $\rho$ and $\mu$ are the fluid's density and dynamic viscosity, respectively.

Now, we will consider a reduced version of the domain $\Omega \subset \Omega_{full}$
(see Figure \ref{fig:duct_mesh-cut}) with the same number of outlets $K$ and also a
reduced wall surface $\Gamma_{wall} \subset \Gamma^{full}_{wall}$. Let us know denote
the outlets of the reduced domain $\Gamma_1,\dots,\Gamma_m$.  In such a reduced
domain, where we are interested in formulating two different reduced models to
approximate the full domain solution as we will discuss in the next sections.

\subsection{The duct boundary condition (DuBC)}

MAPDD assumes that the computational domain consists of arbitrarily
shaped regions -- called \textit{junctions} -- connected by thin ducts. In the
junctions, the full three-dimensional incompressible Navier--Stokes equations are
solved. Within the ducts, the flow velocity (and the test function) is forced to be
parallel to and constant in the axial direction of the duct.

That is, we decompose the full domain as
\begin{equation}
    \Omega_{\text{full}} = \Omega \cup \left( \bigcup_{m=1}^{K}
    \Omega_{\text{duct},m} \right).
\end{equation}
We then write the weak form of the Navier--Stokes equations over this full domain:
\begin{equation}
    \int_{\Omega} \mathcal{F}(\vec u, p, \vec v, q)
    + \sum_{m=1}^K \int_{\Omega_{\text{duct},m}} \mathcal{F}(\vec u, p, \vec v, q)
    = 0,
\end{equation}
with
\begin{equation}
    \mathcal{F}(\vec u, p, \vec v, q) = \rho \frac{\partial \vec{u}}{\partial t}
    \cdot \vec{v}
    + \rho (\vec{u} \cdot \nabla) \vec{u} \cdot \vec{v}
    + \mu \nabla \vec{u} : \nabla \vec{v}
    - p \nabla \cdot \vec{v}
    + q \nabla \cdot \vec{u}.
\end{equation}

In each $\Omega_{\text{duct},m}$, MAPDD enforces the flow assumptions described
above. As a consequence, the pressure, convective and divergence terms vanish.
Taking next the test function also constant along the axial direction of
$\Omega_{\text{duct},m}$, the integral over each duct is reduced to a surface
integral at the interface $\Gamma_m$ between the junction and the duct, times the
duct's length $\ell_m$, as follows:
\begin{equation}
    \sum_{m=1}^K \int_{\Omega_{\text{duct},m}} \mathcal{F}(\vec u, p, \vec v, q)
    =
    \sum_{m=1}^K \ell_m \left\{
        \int_{\Gamma_m} \rho \frac{\partial u_n}{\partial t} v_n
        + \mu \nabla_t u_n \cdot \nabla_t v_n
    \right\},
    \label{eq:mapdd_surface_term}
\end{equation}
where $u_n = \vec{u} \cdot \vec{n}$ and $v_n = \vec{v} \cdot \vec{n}$ are the
normal components of the velocity and test function, respectively. The operator is given by
$\nabla_t(\cdot) = \nabla(\cdot) - (\nabla(\cdot) \cdot \vec{n}) \vec{n}$.


Consequently, the incompressible Navier-Stokes equations plus its boundary
conditions, in weak form (for more detail see Lemma 4 of
\cite{bertoglio2019junction}) can be written as: Find $(\vec u(t) ,p(t)) \in
[H^1(\Omega)]^3 \times L^2(\Omega)$ such that:
\begin{subequations}\label{eq:MAPDD_weak}
    \begin{align}
        \int_{\Omega} \rho \frac{\partial \vec{u}}{\partial t} \cdot \vec{v} +
        \rho \left ( \vec{u}  \cdot \nabla \right ) \vec{u} \cdot \vec{v} + \mu
        \nabla \vec{u} : \nabla \vec{v} - p \nabla \cdot \vec{v} + q \nabla \cdot
        \vec{u} \notag \\
        + \sum_{m=1}^K \frac{\rho}{2} \int_{\Gamma_{m}}  \vert \vec{u} \cdot \vec{n} \vert_-
        (\vec{u}\cdot \vec{v})  + \sum_{m=1}^K \ell_m \ \bigg \{  \int_{\Gamma_m} \rho
            \frac{\partial u}{\partial t} v  + \mu \nabla_{t} u_n \cdot \nabla_{t} v_n
        \bigg \} = 0 \\
        \vec{u} \times \vec{n} = \vec{0} \quad \text{on} \quad \{\Gamma_{1},\dots\Gamma_{m}\}
        \times (0,T] \label{eq:enforce_normal_vel} \\
        \vec{u}  =  \vec{u}_{inlet} \quad \text{on} \quad \Gamma_{inlet} \times (0,T] \\
        \vec{u}  = \vec{0} \quad \text{on} \quad \Gamma_{wall} \times (0,T]
    \end{align}
\end{subequations}

for all $(\vec v,q) \in [H_0^1(\Omega)]^3 \times L^2(\Omega)$.

Apart from the standard backflow term, Equation \eqref{eq:enforce_normal_vel} imposes
the velocity to be perpendicular to the outlet, as it is required in the MAPDD theory
in \cite{bertoglio2019junction}.

For the sake of simplicity and comparison with the full domain solution and 3D-0D
approach, we also discretize Problem \eqref{eq:MAPDD_weak} with a backward Euler time
discretization, obtaining the following problem for $k>0$: 
Find $(\vec{u}^{k+1},p^{k+1}) \in [H^1(\Omega)]^3 \times L^2(\Omega)$ such that:


\begin{subequations}\label{eq:MAPDD_weak_discrete}
    \begin{align}
        \int_{\Omega} \rho \left(\frac{\vec{u}^{k+1} - \vec{u}^{k} }{\tau}  + \vec{u}^{k}
            \cdot \nabla \vec{u}^{k+1} + \frac{1}{2} (\nabla \cdot \vec{u}^{k})\vec{u}^{k+1}
        \right) \cdot \vec{v} +  \int_{\Omega} \mu \nabla \vec{u}^{k + 1} : \nabla \vec{v}  \\
        \int_{\Omega} \left( q \nabla \cdot \vec{u}^{k+1} -  p^{k+1} \nabla \cdot \vec{v}
        \right) + \sum_{m=1}^K \int_{\Gamma_{m}} \frac{\rho}{2} \vert \vec{u}^{k} \cdot
        \vec{n} \vert_- (\vec{u}^{k+1} \cdot \vec{v}) \\
        + \sum_{m=1}^K \int_{\Gamma_{m}} \ell_m \left(  \rho  \frac{ u_n^{k+1} - u_n^{k}  }{\tau}
        v  + \mu \nabla_{t} u_n^{k+1} \cdot \nabla_{t} v_n  \right)
        \\
        + \delta_{stream} \int_{\Omega} \left ( \vec{u}^k \cdot \nabla
        \vec{u}^{k+1}\right ) \cdot \left ( \vec{u}^k \cdot \nabla \vec{v} \right )
        \label{eq:MAPDD_weak_streamline} \\
        + \sum_{m=1}^K\int_{\Gamma_{m}}
        \gamma_{tan}  \left ( \vec{u}^{k+1} - u_n^{k+1} \vec{n} \right )  \cdot  \vec{v} = 0
        \label{eq:MAPDD_weak_discrete_integrals}\\
        \vec{u}^{k+1}   =  \vec{u}_{inlet}(t^{k+1}) \quad \text{on} \quad \Gamma_{inlet}  \\
        \vec{u}^{k+1}   = \vec{0} \quad \text{on} \quad \Gamma_{wall}
    \end{align}
\end{subequations}
for all $(\vec v,q) \in [H_0^1(\Omega)]^3 \times L^2(\Omega)$.
The last term in Equation \eqref{eq:MAPDD_weak_discrete_integrals} enforces the fluid
to flow perpendicular to the outlet. Additionally, in
\eqref{eq:MAPDD_weak_streamline}, we include a streamline diffusion stabilization
term, using the same coefficient as in \cite{garay2022parameter}.

Note that in spite of the fully implicit evaluation of the boundary integral terms
consistently with the discretization of the inertial and viscous terms in $\Omega$ --
only local coupling is introduced, in contrast to implicit 3D-0D coupling models.  As
a consequence, it is straightforward to show the unconditional stability of this
approach. Indeed, by  testing Equation \eqref{eq:MAPDD_weak_discrete} in the unforced
case with $\vec{v} = \vec{u}^{k+1}$, we obtain the energy balance:
\begin{align}
    \frac{E^{k+1} - E^{k}}{\tau} ={} & - \int_{\Omega} \mu \| \nabla \vec{u}^{k+1} \|^2 -
    \int_{\Omega} \frac{\rho}{2\tau} \| \vec{u}^{k+1} - \vec{u}^{k} \|^2 \notag \\
    &  - \sum_{m=1}^K \ell_m \bigg \{  \int_{\Gamma_m} \frac{\rho}{2\tau} \| u_n^{k+1} -
    u_n^{k} \|^2 + \int_{\Gamma_m}  \mu \| \nabla_{t} u_n^{k+1} \|^2  \bigg \}  \notag \\
    & -  \int_{\Gamma_{m}}  \frac{\rho}{2} \vert \vec{u}^{k} \cdot \vec{n} \vert_+ \|
    \vec{u}^{k+1}\|^2  - \gamma_{tan} \|\vec{u}^{k+1} - u_n^{k+1} \vec{n} \|^2
    \label{eq:energy_MAPDD_disc2}
    \\
    & \leq 0
\end{align}
with
\begin{equation}
    E^k = \int_{\Omega} \frac{\rho}{2}  \| \vec{u}^k  \|^2 + \sum_{m=1}^{K} \ell_{m}
    \int_{\Gamma_m} \frac{\rho}{2 } ( u_n^k)^2,
\end{equation}

\begin{remark}
    Though the Duct boundary condition does not consider convection in its derivation
    , the introduction of a backflow stabilization is required in theory and practice
    since still the duct boundary is \textit{open}. However, note that the
    viscous term on $\Gamma_m$ has backflow stabilization properties as it was proven in
    \cite{bertoglio2016stokes}. Therefore, a very large value of $\ell_m$ may suffice to
    avoid backflow instabilities.

\end{remark}

\subsection{Numerical experiments}
\label{section:implementation_details}
The goal of this section is to showcase the stability properties of the DuBC method in
a realistic testcase.




\subsubsection{Setup}

\begin{figure}[!hbtp]
    \centering
    \subfloat[]{
        \includegraphics[trim=50 0 0 0, clip,
    width=0.5\textwidth]{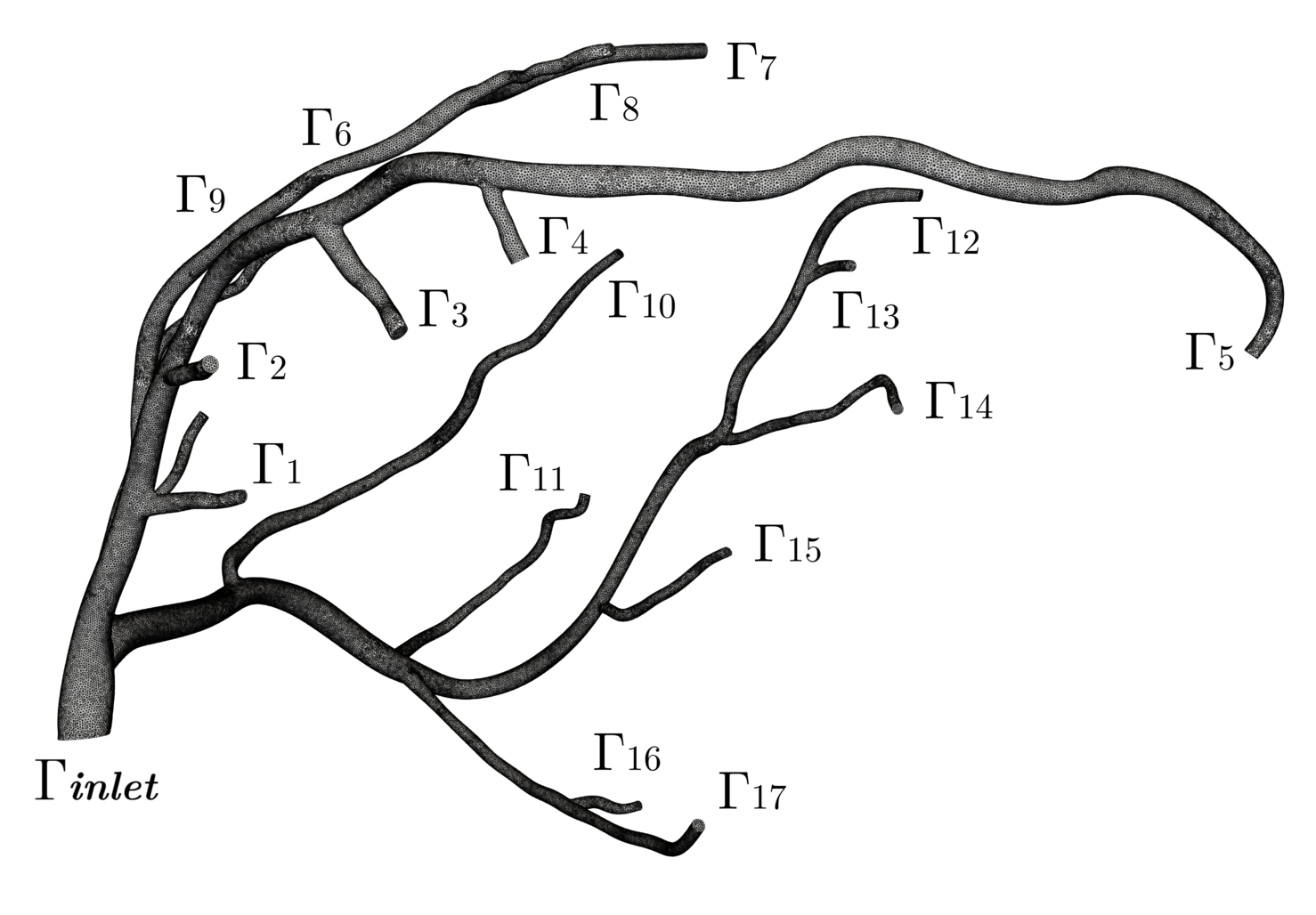} \label{fig:coronary_mesh-full}}
    \subfloat[]{
        \hspace{0.2cm}
        \includegraphics[trim=50 0 0 0, clip,
        width=0.4\textwidth]{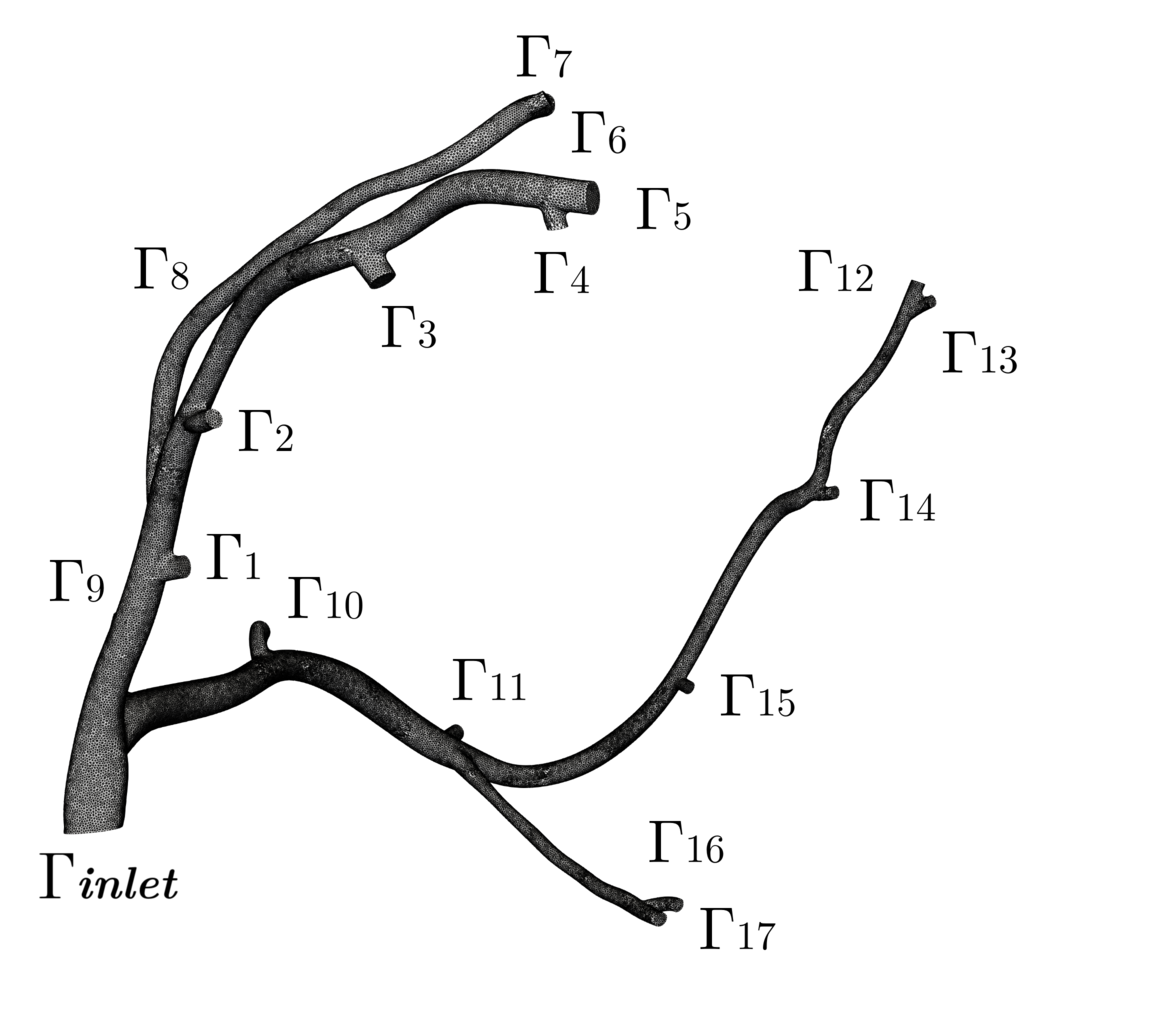}
    \label{fig:coronary_mesh-reduced} }
    \caption{(a) Left coronary artery's full domain $\Omega_{full}$ used for this
        study. (b) Reduced geometry $\Omega$ obtained after cutting the branches
    perpendicularly to their centerline}
    \label{fig:coronary_mesh}
\end{figure}

\paragraph{Geometry and physical constants} We assume a Newtonian fluid with constant
density and dynamic viscosity as $\rho = 1.06 \ \rm{g/cm}^3$ and $\mu = 0.035$ P,
respectively. For the geometry, we consider a left coronary tree with $K=17$ outlets,
see Figure \ref{fig:coronary_mesh-full}.


Next, a reduced model was obtained from the original coronary geometry by cutting the
segments perpendicular to their centerline where all branches do not lead to a
bifurcation, see Figure \ref{fig:coronary_mesh-reduced}. Consequently, the obtained
reduced geometry also possess the same amount of outlets than the original, but with a
total volume reduction of around $66 \%$. The approximate lengths were used as
parameters for the DuBC, representing the missing duct-like structures that we neglect
and simulate through the model itself. The resulting values are summarized in Table
\ref{tab:Parameters}.



\begin{table}[!hbtp] \centering
    \footnotesize
    \sisetup{table-format=1, table-number-alignment=center}
    \begin{tabular}{cS[table-format=1]*{9}{S}}
        \toprule
        \addlinespace
        {Boundary} &  & {$\Gamma_1$} & {$\Gamma_2$} & {$\Gamma_3$} & {$\Gamma_4$} &
        {$\Gamma_5$} & {$\Gamma_6$} & {$\Gamma_7$} & {$\Gamma_8$} & {$\Gamma_9$}
        \tabularnewline
        \cmidrule[\lightrulewidth](lr){1-11}\addlinespace[1ex]
        $\ell_m \ \rm{(cm)}$ & & {0.95} & {0.47}  & {0.98} & {0.46}  & {7.94}  &  {1.38} &
        {2.80} & {4.85} & {2.39}  \tabularnewline \addlinespace[1ex]
        \cmidrule[\lightrulewidth](lr){1-11}\addlinespace[1ex]
        {Boundary} &  & {$\Gamma_{10}$} & {$\Gamma_{11}$} & {$\Gamma_{12}$} &
        {$\Gamma_{13}$} & {$\Gamma_{14}$} & {$\Gamma_{15}$} & {$\Gamma_{16}$} & {$\Gamma_{17}$}
        \tabularnewline
        \cmidrule[\lightrulewidth](lr){1-11}\addlinespace[1ex]
        $\ell_m \ \rm{(cm)}$ & & {6.84} & {3.21} & {1.63}  & {0.42}  & {2.70}  & {1.99} &
        {0.51} & {2.48} \tabularnewline \addlinespace[1ex]
        \bottomrule
    \end{tabular}
    \caption{Virtual lengths $\ell$ parameters in centimeters of the MAPDD model, for
    every open boundary in the coronary arteries.}
    \label{tab:Parameters}
\end{table}

\paragraph{Temporal discretization} 

Three different values were chosen for the time-step as $\tau = 0.001, \ 0.005$ and
$0.01 s$, with a total simulation time of $0.9 \ s$. The initial conditions were set
as $\vec{u}^0 =\vec{0}$.


\paragraph{Spatial discretization} The weak forms of the incompressible Navier-Stokes
equations where discretized using stabilized $\mathbb{P}1/\mathbb{P}1$ Taylor-Hood
elements. The full computational mesh consisted in 794,705 tetrahedrons and 187,566
vertices. Once reduced, the resulting mesh had a total of 508,970 tetrahedrons and
113,637 vertices.



\paragraph{Dirichlet boundary conditions via penalization}
The inlet Dirichlet boundary condition is introduced in the variational form as a
penalization term $A_{inlet}$ defined as:
\begin{equation}
    A_{inlet} = \gamma_{inlet} \int_{\Gamma_{inlet}} (\vec{u} - \vec{u}_{inlet}) \cdot \vec{v}
\end{equation}

where the parameter $\gamma_{inlet}$ was fixed in $10^5 \ g/(cm^2 \cdot s)$, while
$\vec{u}_{inlet}$ is set as:
\begin{equation}
    \vec{u}_{inlet}(\vec{x},t) = f(t) \ \vec{u}_{stokes}(\vec{x}),
\end{equation}
where $\vec{u}_{stokes}(\vec{x})$ is the solution of a steady Stokes problem inside
the domain, having a parabolic-like shape adapted to the mesh geometry. The function
$f(t)$ stands for the time-dependency of the inflow velocity and is taken in such a
way that the total flow through $\Gamma_{inlet}$ follows a population-averaged curve
taken from \cite{pappano201310}.

Finally, the tangential penalization term introduced in Equation
\eqref{eq:MAPDD_weak_discrete} $\gamma_{tan}$ was fixed in $10^8 \ \rm{g}/(\rm{cm}^2
\cdot \rm{s})$ for all cases.

\subsubsection{Results}

\begin{figure}[htbp]
    \centering
    \makebox[.3\textwidth]{}\hfill
    \subfloat[ Reference at $\tau = 1 \rm{ms}$]{
        \includegraphics[trim=100 5 175 25, clip,
    width=0.33\textwidth]{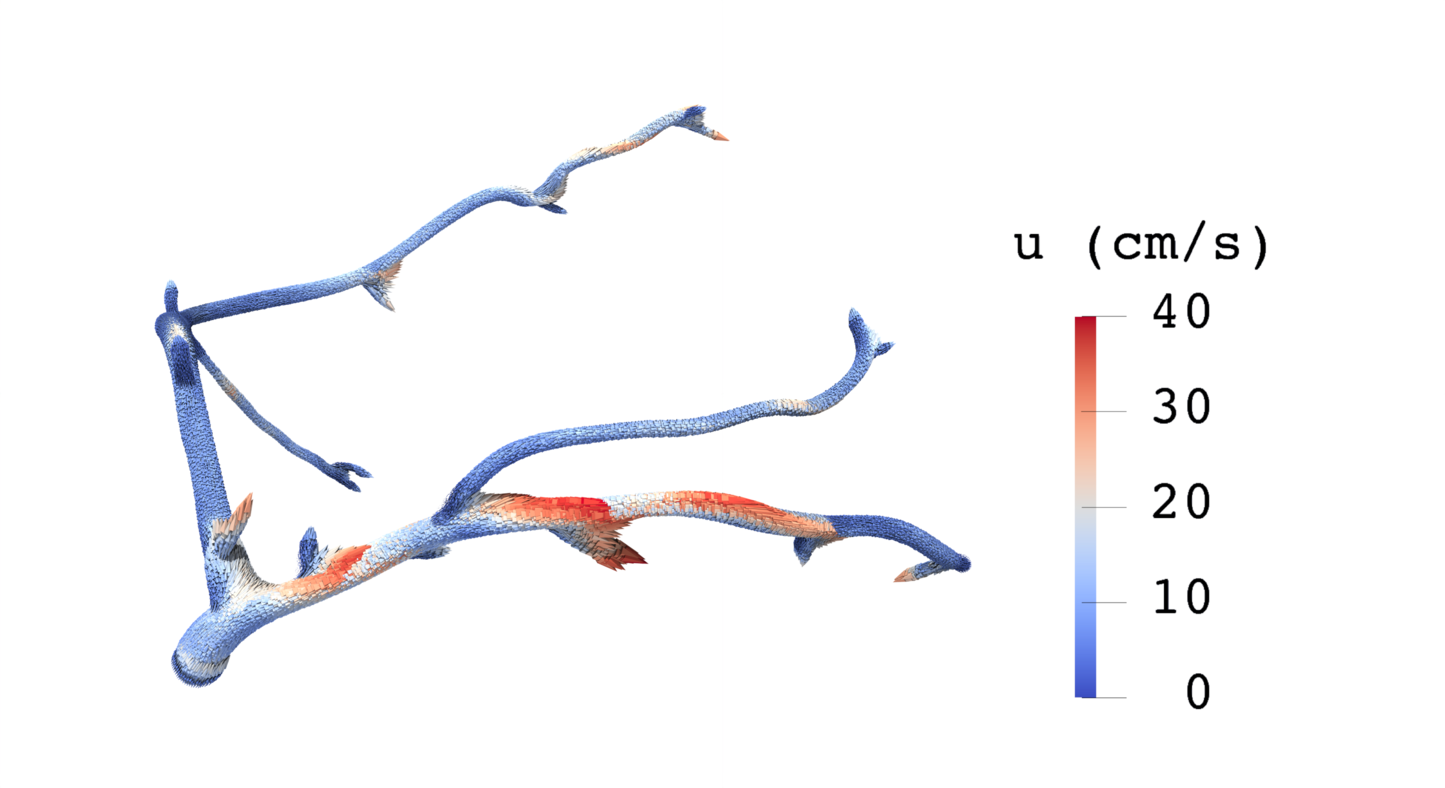} }
    \makebox[.3\textwidth]{}\hfill
    \subfloat[DuBC at $\tau = 1 \rm{ms}$]{
        \includegraphics[trim=100 5 175 25, clip,
    width=0.33\textwidth]{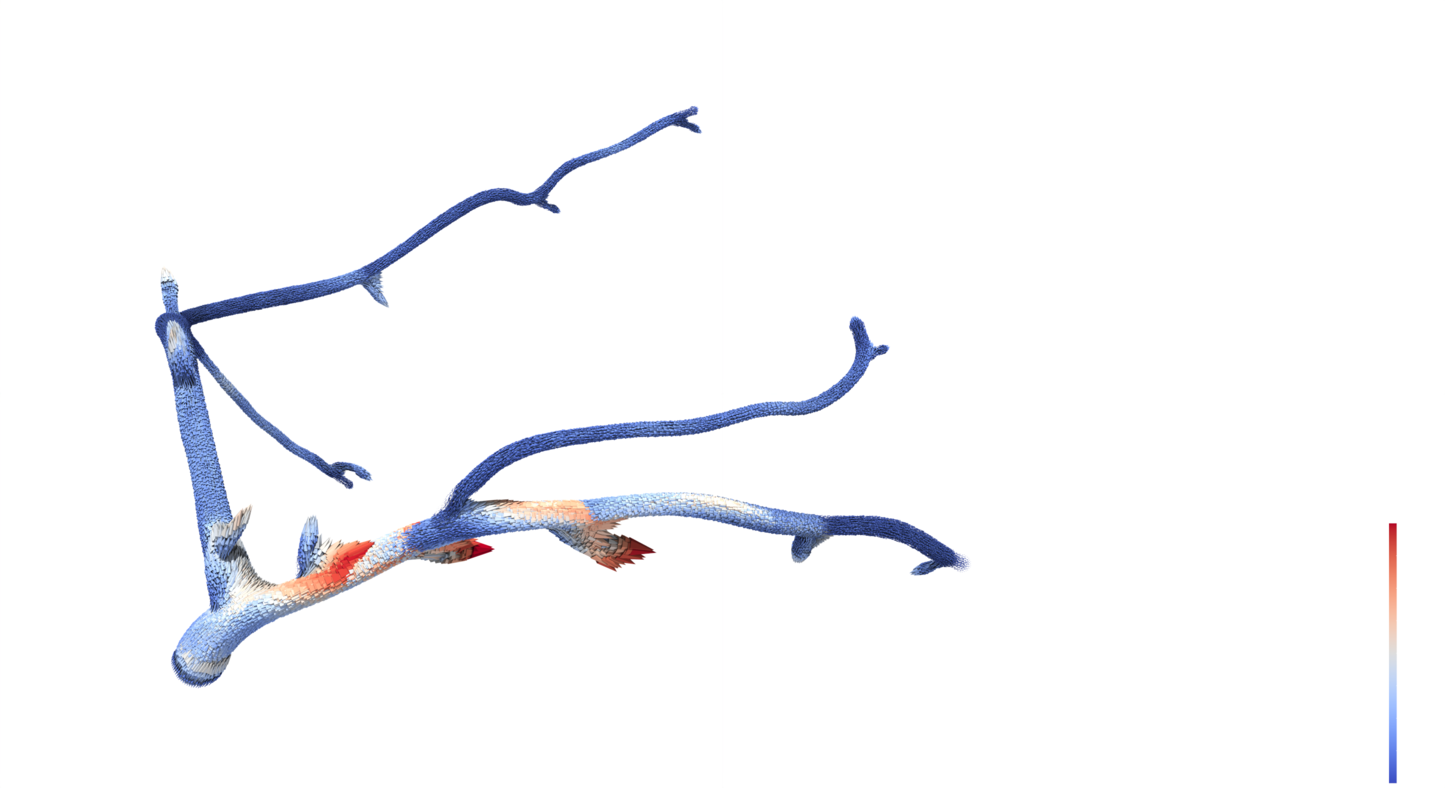} }
    \subfloat[DuBC at $\tau = 5 \rm{ms}$]{
        \includegraphics[trim=100 5 175 25, clip,
    width=0.33\textwidth]{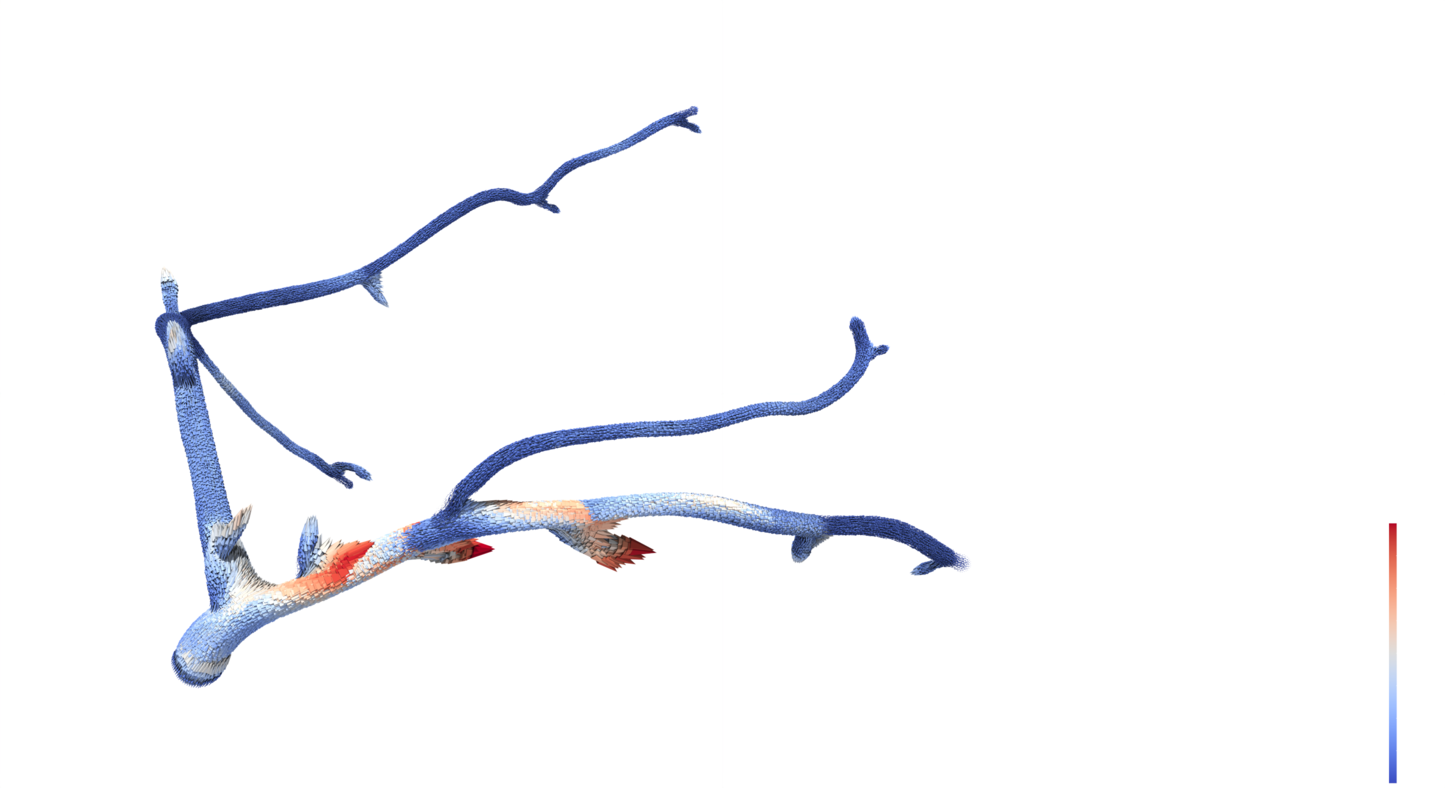} }
    \subfloat[DuBC at $\tau = 10 \rm{ms}$]{
        \includegraphics[trim=100 5 175 25, clip,
    width=0.33\textwidth]{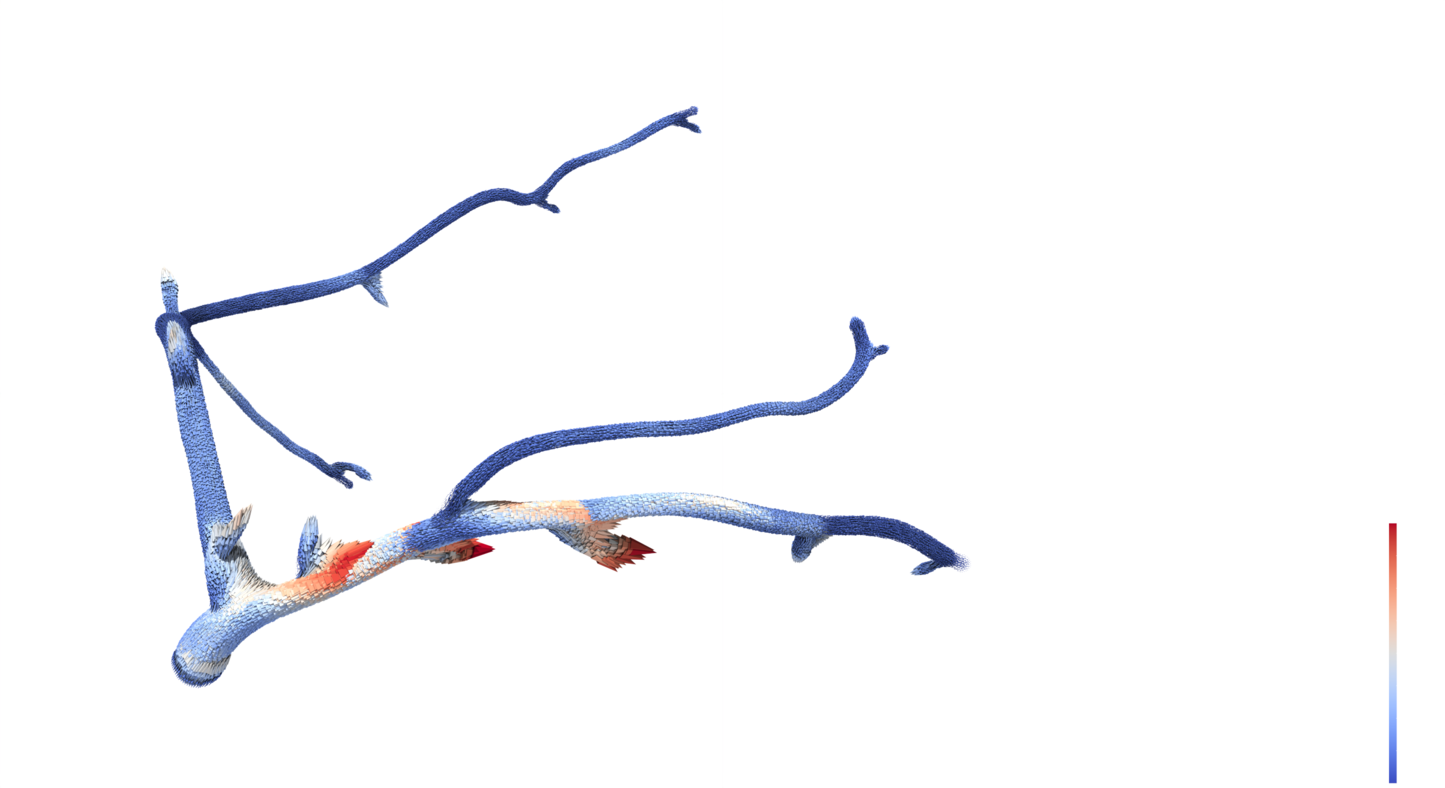} }\\
    \caption{Velocity fields when using different models with increasing simulation time
    step. All colormaps and arrow sizes are the same}
    \label{fig:velocities_at_different_dt}
\end{figure}

We first compare the full domain solution of Section \ref{sec:reference_solution} on
the entire geometry against the solutions on the reduced geometry using the DuBC. In
order to do that, we interpolate the full domain solution onto the reduced geometry to
make all solutions comparable. Figure \ref{fig:velocities_at_different_dt} shows the
velocities obtained at peak ($t = 0.69 \rm{s}$) for the DuBC model when
using different simulation's time steps. From these, it can be observed that the
velocities with  DuBC are highly robust to increases in the time step.

\begin{figure}[hbtp]
    \centering
    \includegraphics[trim=0 0 0 0, clip,
    width=0.5\textwidth]{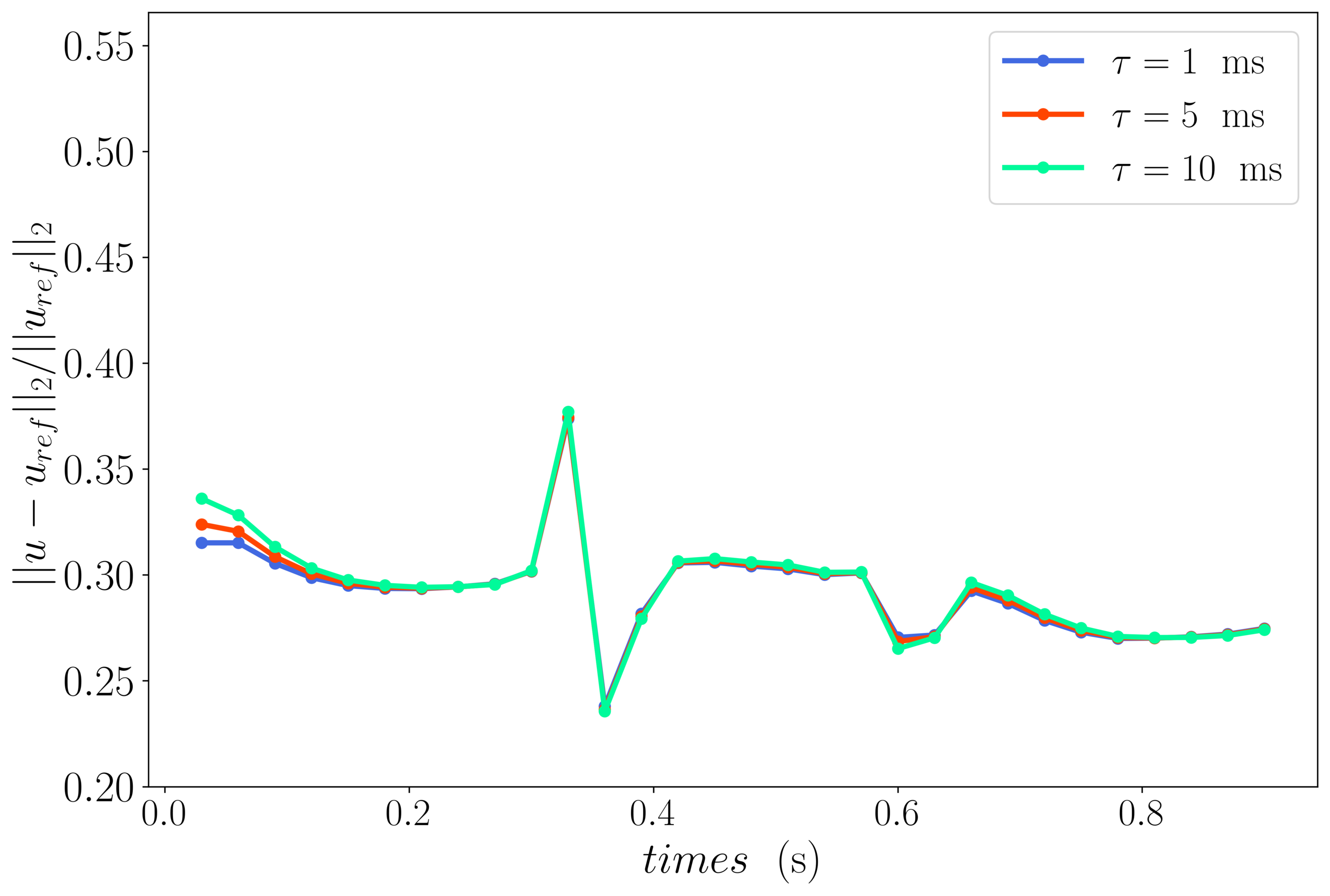}
    \caption{Comparison of the relative $L_2$ norms of the velocity obtained with the
    DuBC model as the simulation time step increases.}
    \label{fig:norms_by_model}
\end{figure}

Additionally, Figure \ref{fig:norms_by_model} shows the $L_2$ norms of the
velocities relative to a reference solution computed as:
\begin{equation}
    \epsilon(t) = \displaystyle \frac{\sum_i | \vec{u}(t) - \vec{u}_{ref}(t)
    |^2}{\sum_i | \vec{u}_{ref}(t)  |^2}
    \label{ec:norms_formula}
\end{equation}
where the summation is over all nodes of the mesh. For simplicity, the error was only
computed every $0.03$s. The reference is obtained solving the full domain with $\tau =
1$ ms and interpolating to the reduced geometry afterwards. From these curves, it can
be observed that the impact of increasing the simulation time-step 10 fold produced
little impact of the
solution quality.

\section{Fractional step discretization including DuBC} \label{sec:CT_algo}
\subsection{Formulation}

The MAPDD method was originally developed and numerically tested on monolithic schemes
\cite{panasenko1998method,blanc1999asymptotic,bertoglio2019junction,bertoglio2021reconstruction},
leading to the
framework to solve the DuBC as proposed and analyzed in the previous section.

Now, we will extend the DuBC formulation to a non-incremental fractional step
discretization in order to accelerate the computational time. We start from a version
of the classical Chorin-Temam non-incremental pressure correction scheme
\cite{guermond2006overview}.

In fractional step approaches, velocity and pressure solutions are staggered, in the
so-called \textit{tentative} and \textit{(pressure) projection} steps. In the case of
the Chorin-Temam method, the former corresponds to the same formulation of the
monolithic problem, with the exception that the pressure is evaluated explicitly.
Therefore, the tentative step when using de DuBC becomes the same as in Equation
\eqref{eq:MAPDD_weak_discrete} but with the $p^k$ instead $p^{k+1}$. Note that such
formulation already considers the \textit{couple first, then discretize} strategy,
with ``discretize" here meaning both spatially and temporally.


However, the pressure projection step requires more attention.
Here, the left-hand-side of the pressure projection step can be rewritten as:
\begin{equation}
    \int_{\Omega_{full}} \nabla p^{k} \cdot \nabla  q = \int_{\Omega} \nabla p^{k} \cdot
    \nabla q  + \sum_{m=1}^K \int_{\Omega_{\text{duct},m}} \nabla p^{k} \cdot \nabla q
    \label{eq:pressure_gradient_raw}
\end{equation}

In the DuBC approach, the fluid pressure in the extension $\Omega_{\text{duct},m}$
is assumed to be a constant gradient along the duct, where the pressure at the end is
zero. Consequently, the last term of Equation \eqref{eq:pressure_gradient_raw} can be
rewritten as:

\begin{equation} \label{eq:pressure_gradient_step01}
    \sum_{m=1}^K \int_{\Omega_{\text{duct},m}} \nabla p^{k} \cdot \nabla q 
    = \sum_{k=0}^{K} \ell_m \int_{\Gamma_m} \frac{-p^{k}}{\ell_m}   \frac{-q}{\ell_m} =
    \sum_{k=0}^{K} \frac{1}{\ell_m} \int_{\Gamma_m}  p^{k} q
\end{equation}
becoming a penalization of the pressure according to the length of the extension.

The complete algorithm is detailed in Algorithm \ref{alg:mapdd_CT}. Note that at the
projection step, there is no integral of the divergence of the velocity in
$\Omega_{\text{duct},m}$, since by construction of the MAPDD (prior to temporal discretization),
the velocity trial and test functions are divergence free in the extension.

\begin{algorithm}[H]
    \caption{Fractional step algorithm with DuBC}
    \label{alg:mapdd_CT}
    Given $\vec{u}^0 \in [H^1(\Omega)]^3$, perform for $n\geq0$: \\

    \textbf{1. Pressure projection step:}
    Find $p^{k} \in H^1(\Omega)$ such that:
    \begin{align}\label{eq:pressurestep}
        \int_\Omega \nabla p^{k} \cdot \nabla q  + \frac{\rho}{\tau}  \int_\Omega \nabla \cdot
        \vec{u}^{k} q + \sum_{k =1}^\nbOut  \int_{\Gamma_m}  \frac{p^{k} q}{\ell_m}   = 0
    \end{align}
    for all $q \in H^1(\Omega)$.

    \textbf{2. Tentative 
    velocity step:} Find $\vec{u}^{k+1} \in
    [H^1(\Omega)]^3$ such that:
    \begin{subequations}
        \begin{align}
            \int_{\Omega} \rho \left(\frac{\vec{u}^{k+1} - \vec{u}^{k} }{\tau}  + \vec{u}^{k}
                \cdot \nabla \vec{u}^{k+1} + \frac{1}{2} (\nabla \cdot \vec{u}^{k})\vec{u}^{k+1}
            \right) \cdot \vec{v} + \int_{\Omega} \mu \nabla \vec{u}^{k + 1} : \nabla \vec{v}  -
            p^{k} \nabla \cdot \vec{v} \notag \\
            + \sum_{m=1}^K \int_{\Gamma_{m}} \frac{\rho}{2} \vert \vec{u}^{k} \cdot \vec{n}
            \vert_- (\vec{u}^{k+1} \cdot \vec{v}) + \ell_m \left(  \rho  \frac{
                    u_n^{k+1} - u_n^{k}
            }{\tau} v  + \mu \nabla_{t} u_n^{k+1} \cdot \nabla_{t} v_n  \right)
            \notag \\
            + \sum_{m=1}^K\int_{\Gamma_{m}}
            \gamma_{tan}  \left ( \vec{u}^{k+1} - u_n^{k+1} \vec{n} \right )  \cdot  \vec{v} = 0
            \label{eq:MAPDD_weak_discrete_integrals_CT}\\
            \vec{u}^{k+1}   =  \vec{u}_{inlet}(t^{k+1}) \quad \text{on} \quad \Gamma_{inlet}  \\
            \vec{u}^{k+1}   = \vec{0} \quad \text{on} \quad \Gamma_{wall}
        \end{align}
    \end{subequations} \label{eq:viscousstep}
    for all $\mathbf{v} \in [H^1_0(\Omega)]^3$.

\end{algorithm}



The energy balance in the non-forced case can be obtained by testing with
$\vec{v}=\vec{u}^{k+1}$ as well as $q=p^{k}$:
\begin{align}
    \frac{E^{k+1} - E^{k}}{\tau} = & - \int_{\Omega} \mu \| \nabla \vec{u}^{k+1} \|^2 -
    \int_{\Omega} \frac{\rho}{2\tau} \| \vec{u}^{k+1} - \vec{u}^{k} \|^2 + \int_\Omega p^k
    \nabla \cdot \vec{u}^{k+1}  \notag \\
    &  - \sum_{m=1}^K \ell_m \bigg \{  \int_{\Gamma_m} \frac{\rho}{2\tau} \| u_n^{k+1} -
    u_n^{k} \|^2 + \int_{\Gamma_m}  \mu \| \nabla_{t} u_n^{k+1} \|^2  \bigg \}  \notag \\
    & -  \int_{\Gamma_{m}}  \frac{\rho}{2} \vert \vec{u}^{k} \cdot \vec{n} \vert_+ \|
    \vec{u}^{k+1}\|^2  - \gamma_{tan} \|\vec{u}^{k+1} - u_n^{k+1} \vec{n} \|^2 \notag \\
    & - \int_{\Omega} \frac{\tau}{\rho}  \| \nabla p^{k} \|^{2} - \int_\Omega p^k \nabla
    \cdot \vec{u}^k \notag \\
    & - \sum_{m=1}^K  \frac{\tau}{\rho}\int_{\Gamma_m} \frac{1}{\ell_m}  (p^{k})^2 +
    \gamma_{press} \| \nabla_{t} p^{k} \|^2   \label{eq:energy_fractional_step}
\end{align}
As in standard energy balance proofs for CT methods
\cite{bertoglio2013fractional,guermond2006overview,fernandez2007projection}, we treat
the unsigned terms as follows:
\begin{align}
    \int_\Omega p^k \nabla \cdot (\vec{u}^{k+1} - \vec{u}^k) = & - \int_\Omega \nabla p^k
    \cdot (\vec{u}^{k+1} - \vec{u}^k) + \sum_{m=1}^K \int_{\Gamma_m}  p^k  (u^{k+1} - u^k) \\
    = & \int_\Omega \sqrt{\frac{\tau}{\rho}} \nabla p^k \cdot \sqrt{\frac{\rho}{\tau}}
    (\vec{u}^{k+1} - \vec{u}^k)  + \sum_{m=1}^K \int_{\Gamma_m}
    \sqrt{\frac{\tau}{\rho\ell_m}} p^k  \sqrt{\frac{\ell_m\rho}{\tau}} (u^{k+1} - u^k) \\
    \leq & \int_\Omega \frac{\tau}{2\rho} \| \nabla p^k\|^2 + \frac{\rho}{2\tau}
    \|\vec{u}^{k+1} - \vec{u}^k\|^2  \notag \\
    & + \sum_{m=1}^K \int_{\Gamma_m}  \frac{\tau}{2\rho\ell_m} (p^k)^2  +
    \frac{\ell_m\rho}{2\tau} (u^{k+1} - u^k)^2
\end{align}
which combining with Equation \eqref{eq:energy_fractional_step} results in:
\begin{align}
    \frac{E^{k+1} - E^{k}}{\tau} \leq & - \int_{\Omega} \mu \| \nabla \vec{u}^{k+1} \|^2
    - \sum_{m=1}^K \ell_m \int_{\Gamma_m}  \mu \| \nabla_{t} u_n^{k+1} \|^2  \notag \\
    & -  \int_{\Gamma_{m}}  \frac{\rho}{2} \vert \vec{u}^{k} \cdot \vec{n} \vert_+ \|
    \vec{u}^{k+1}\|^2  - \gamma_{tan} \|\vec{u}^{k+1} - u_n^{k+1} \vec{n} \|^2 \notag \\
    & - \int_{\Omega} \frac{\tau}{2\rho}  \| \nabla p^{k} \|^{2} - \sum_{m=1}^K
    \frac{\tau}{\rho}\int_{\Gamma_m} \frac{1}{2\ell_m}  (p^{k})^2 + \gamma_{press} \|
    \nabla_{t} p^{k} \|^2   \label{eq:energy_fractional_step_final}
\end{align}
leading therefore to unconditional stability of the CT-DuBC formulation.

\subsection{Numerical experiments}

\paragraph{Spatial discretization and details} Algorithm \ref{alg:mapdd_CT} was solved in
the reduced coronary model. As in the monolithic case, a streamline diffusion
stabilization was added and the inflow boundary condition was applied using
penalization. 

\begin{figure}[htbp]
    \centering
    \subfloat[ Monolithic at $\tau = 1 \rm{ms}$]{
    \includegraphics[trim=100 25 125 0, clip, width=0.3\textwidth]{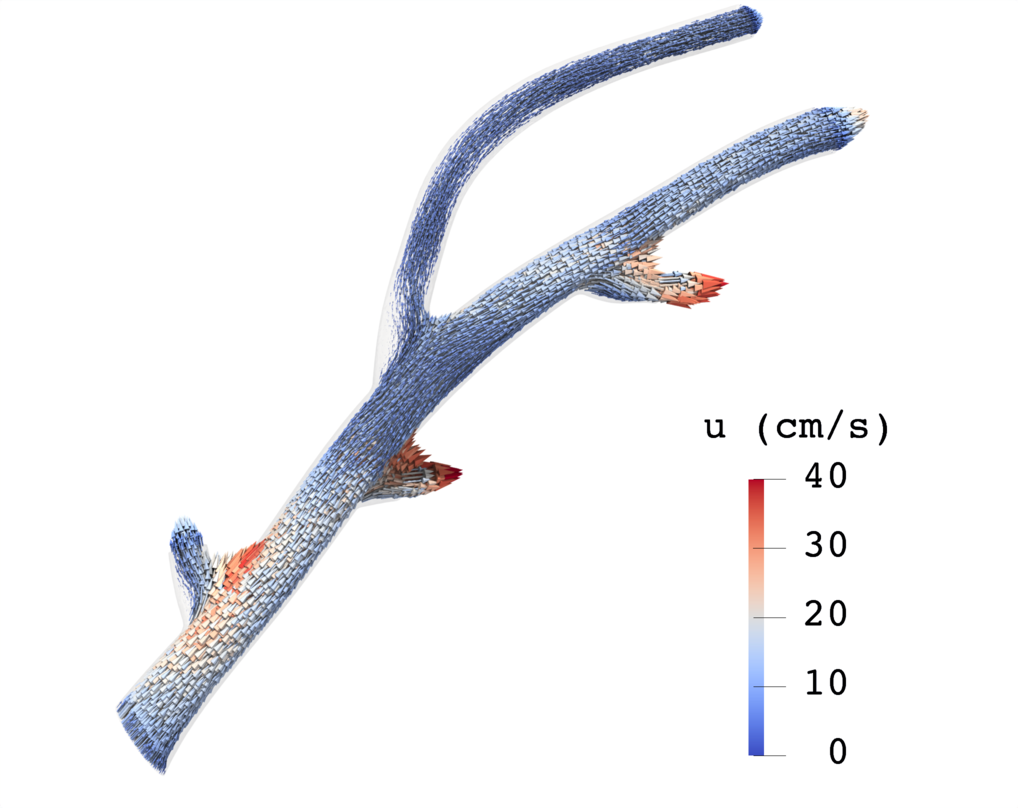} }
    \subfloat[ Monolithic at $\tau = 5 \rm{ms}$]{
        \includegraphics[trim=100 25 125 0, clip,
    width=0.3\textwidth]{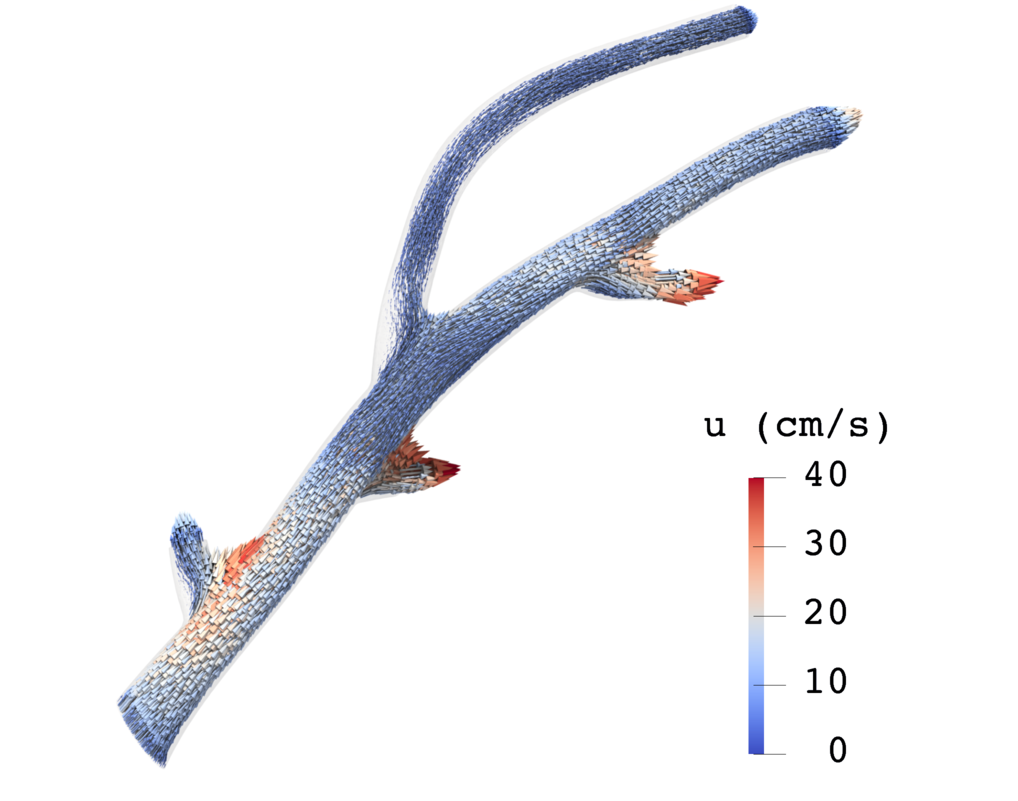} }
    \subfloat[ Monolithic at $\tau = 10 \rm{ms}$ ]{
        \includegraphics[trim=100 25 125 0, clip,
    width=0.3\textwidth]{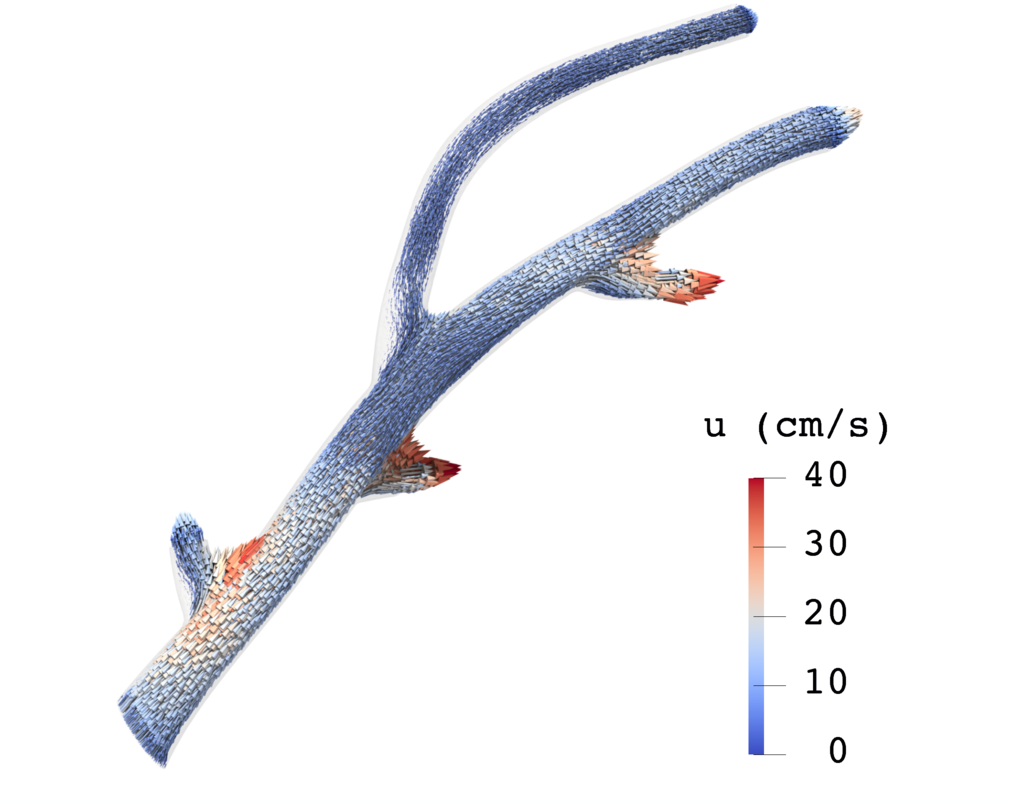} }
    \\
    \subfloat[ CT at $\tau = 1 \rm{ms}$]{
    \includegraphics[trim=100 25 125 0, clip, width=0.3\textwidth]{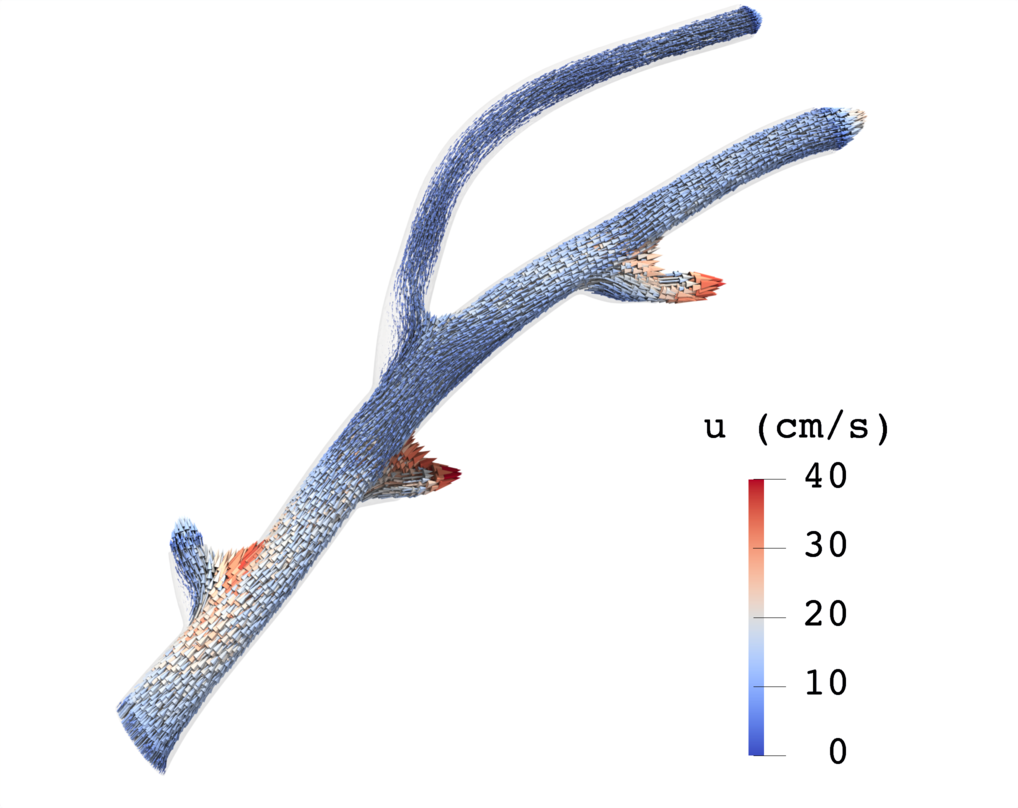} }
    \subfloat[ CT at $\tau = 5 \rm{ms}$]{
        \includegraphics[trim=100 25 125 0, clip,
    width=0.3\textwidth]{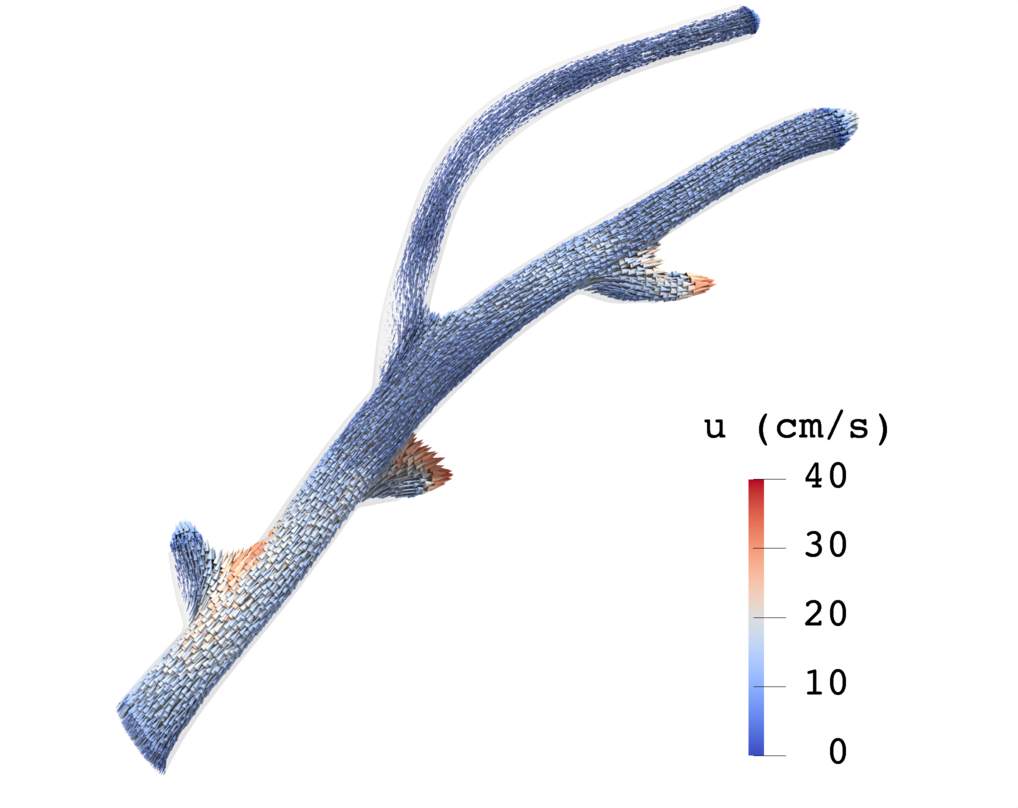} }
    \subfloat[ CT at $\tau = 10 \rm{ms}$ ]{
        \includegraphics[trim=100 25 125 0, clip,
    width=0.3\textwidth]{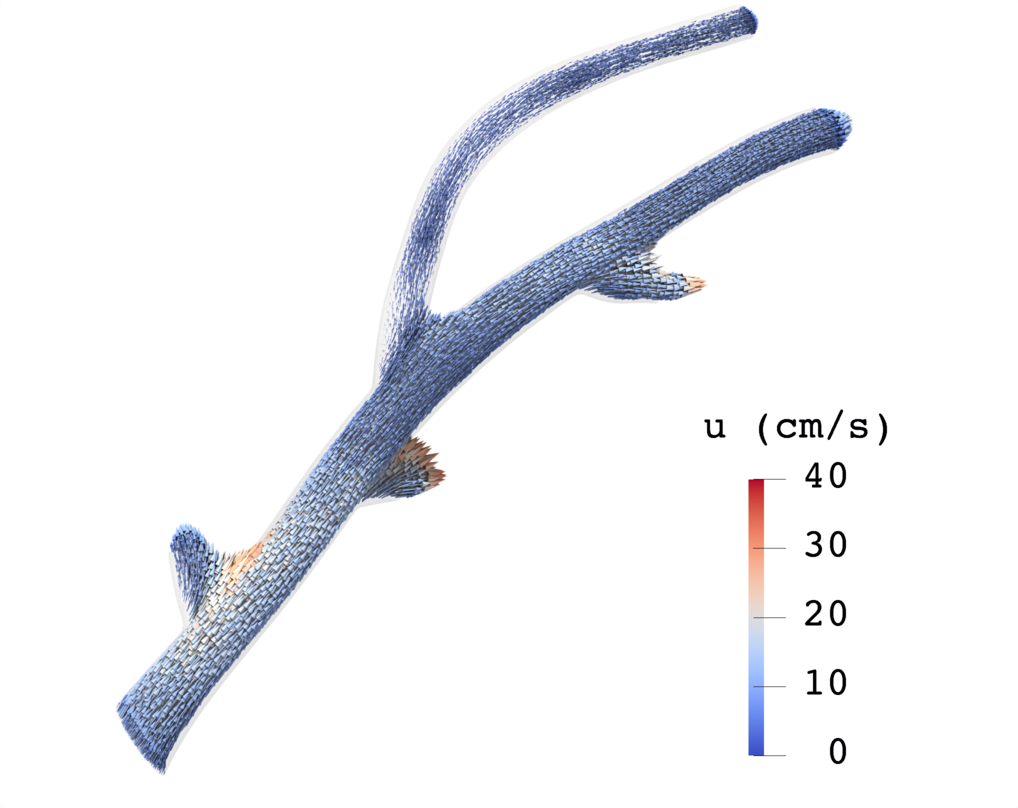} }
    \caption{Velocity fields in the LAD portion of the left coronary artery when
    using the DuBC with monolithic and CT approaches for three different simulation time steps.}
    \label{fig:comparison_ct_dts}
\end{figure}

\paragraph{Results} Figure \ref{fig:comparison_ct_dts} shows, for different simulation
time steps $\tau$, the solutions obtained at the instant of peak inlet velocity ($t =
0.69 \rm{s}$), in a section of the \emph{left anterior descendent} artery (LAD), where
the maximum velocities and also the higher discrepancies against the monolithically
solved DuBC-problem are found. 

\begin{figure}[hbtp]
    \centering
    \includegraphics[trim=0 0 0 0, clip,
    width=0.5\textwidth]{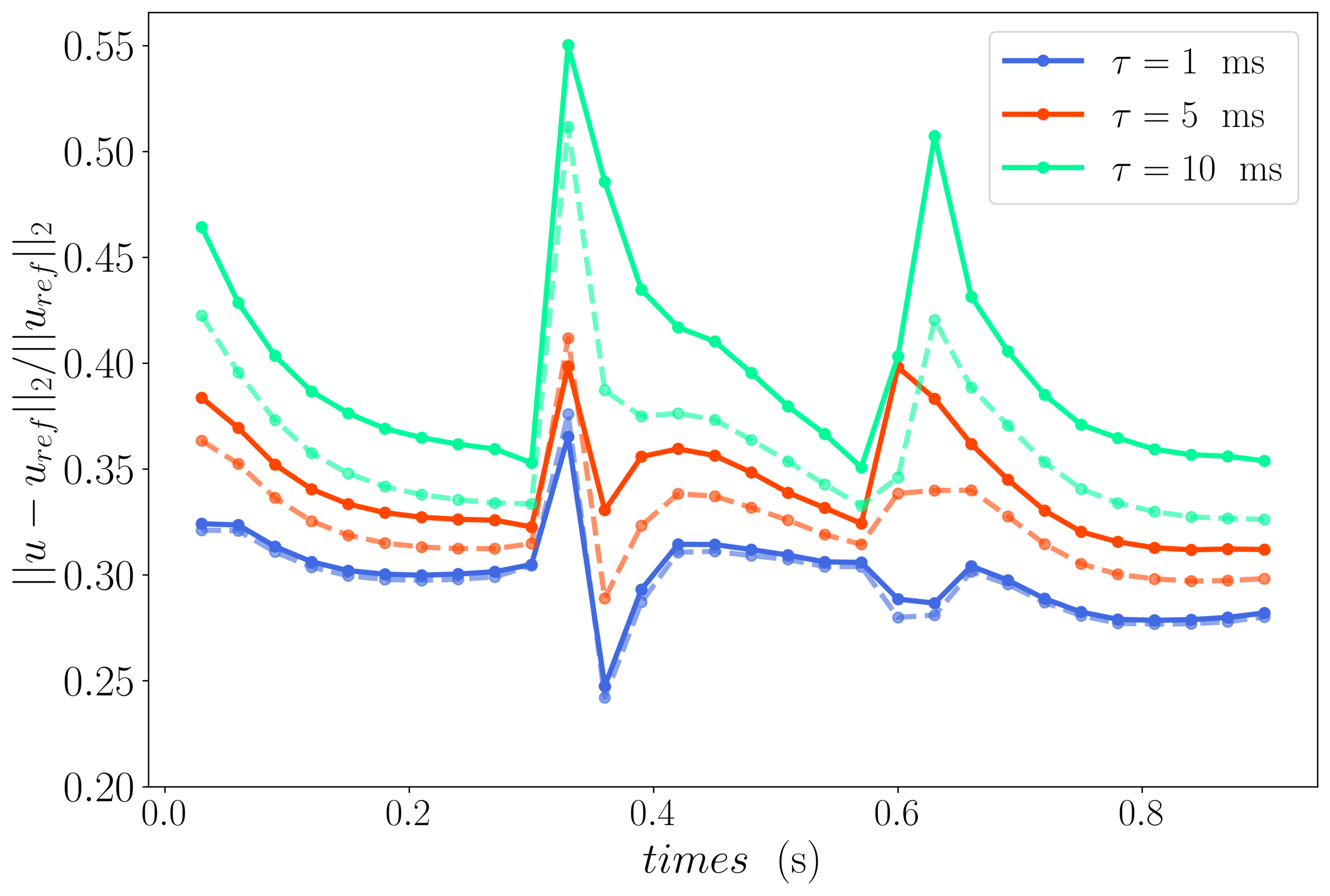}
    \caption{Comparison of the relative $L_2$ norms of the velocity obtained with the
        DuBC model as the simulation time step increases, for the CT method. Dashed lines
        correspond to the the corrected velocity, while continuous lines to the one computed
    using Algorithm \ref{alg:mapdd_CT}.}
    \label{fig:norms_by_model_CT}
\end{figure}

Figure \ref{fig:norms_by_model_CT} shows the relative $L^2$ error norms with respect
to the monolithic reference solution (for both CT velocities) as the time step
increases, computed as in Equation \eqref{ec:norms_formula} at the same time intervals
(every $0.03$s). The results
indicate that the solution remains very close to the monolithic one for smaller values
of $\tau$, but deteriorates as $\tau$ grows. We also include the
error norms associated with the corrected velocity, which is computed by solving an
$L^2$-projection of $\vec{u}^k-\frac{\tau}{\rho}\nabla p^k$, and consistently exhibits
lower errors compared to the main CT velocity.

Finally, Table \ref{tab:forward_runningtimes} shows overall time-steps used and the
resulting running times of all forward simulations, when using 2 cores on a AMD Ryzen
9 7950X with 64 GB RAM. In all cases, a direct LU method was
used for solving the discretized problem. 

\begin{table}[htbp]
    \centering

    \begin{tabular}{c|c|c|c}
        {time step $\tau$} &  {Full model}  &  {DuBC} & {CT-DuBC}\\
        \midrule
        {1ms}  &  {3h\ 4m}  & {2h\ 3m} & {1h\ 11m}\tabularnewline
        \midrule
        {5ms}  &  {34m}  & {24m} & {17m}\tabularnewline
        \midrule
        {10ms}  &  {17m}  & {12m} & {7m} 
    \end{tabular}
    \caption{Total running times for the DuBC model when varying the simulation time-step.}
    \label{tab:forward_runningtimes}
\end{table}



\section{Estimation of DuBC parameters from velocity data} \label{sec:inv}

In this section we present a parameter optimization problem involving the DuBC in an
example of relevance in computational hemodynamics, namely to estimate the lengths
$\ell_1,\dots,\ell_m$ from velocity measurements and the vessel geometry, as they
would be obtained from \textit{4D Flow MRI} \cite{markl20124d,bissell20234d}. The
purpose is to show how the duct boundary condition is well suited for
patient-specific modeling, both in terms of number of parameters to be estimated as
well as its robustness with respect to the parameter values - both crucial features in
parameter estimation problems.


As a parameter estimation method, we employ a Reduced-order Unscented Kalman Filter
(ROUKF) \cite{moireau2011reduced}, which is of wide use in blood flow problems
\cite{bertoglio-moireau-gerbeau-11,moireau-etal-externaltissueestimation-11,bertoglio2014jbmech,nolte_reducing_2019,arthurs2020flexible,garay2022parameter}
and present a computationally tractable way to deal with large time dependent PDE
models as the one used here.

\subsection{Measurement generation}
\label{sec:measurement_gen}

We first define a \emph{high-fidelity} dataset as the solution obtained with the CT
method described in the previous section, at $\tau = 1$ ms. Since our parameter
estimation framework assumes the presence of noise in the measurements, this
dataset was perturbed by adding Gaussian noise with zero mean and a standard
deviation of approximately $5\%$ of the maximum velocity.

\begin{figure}[!hbtp]
    \centering
    \subfloat[Voxel-like mesh used for spatial interpolation]{
        \includegraphics[trim=250 200 150 0, clip,
    width=0.5\textwidth]{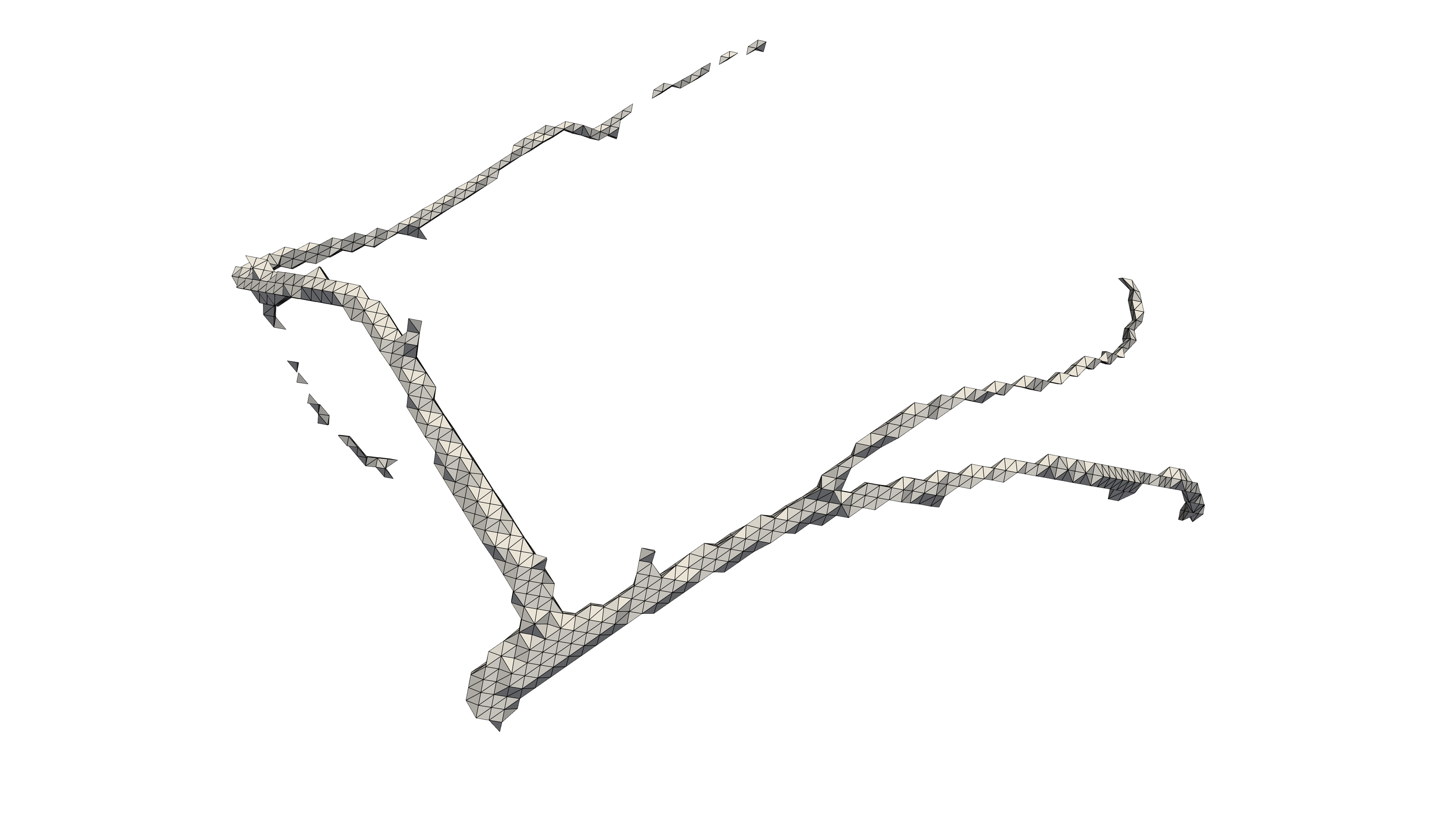}  \label{fig:4d_flow_mesh}}
    \subfloat[Simulated velocity measurements at peak]{
        \includegraphics[trim=250 200 150 0, clip,
    width=0.5\textwidth]{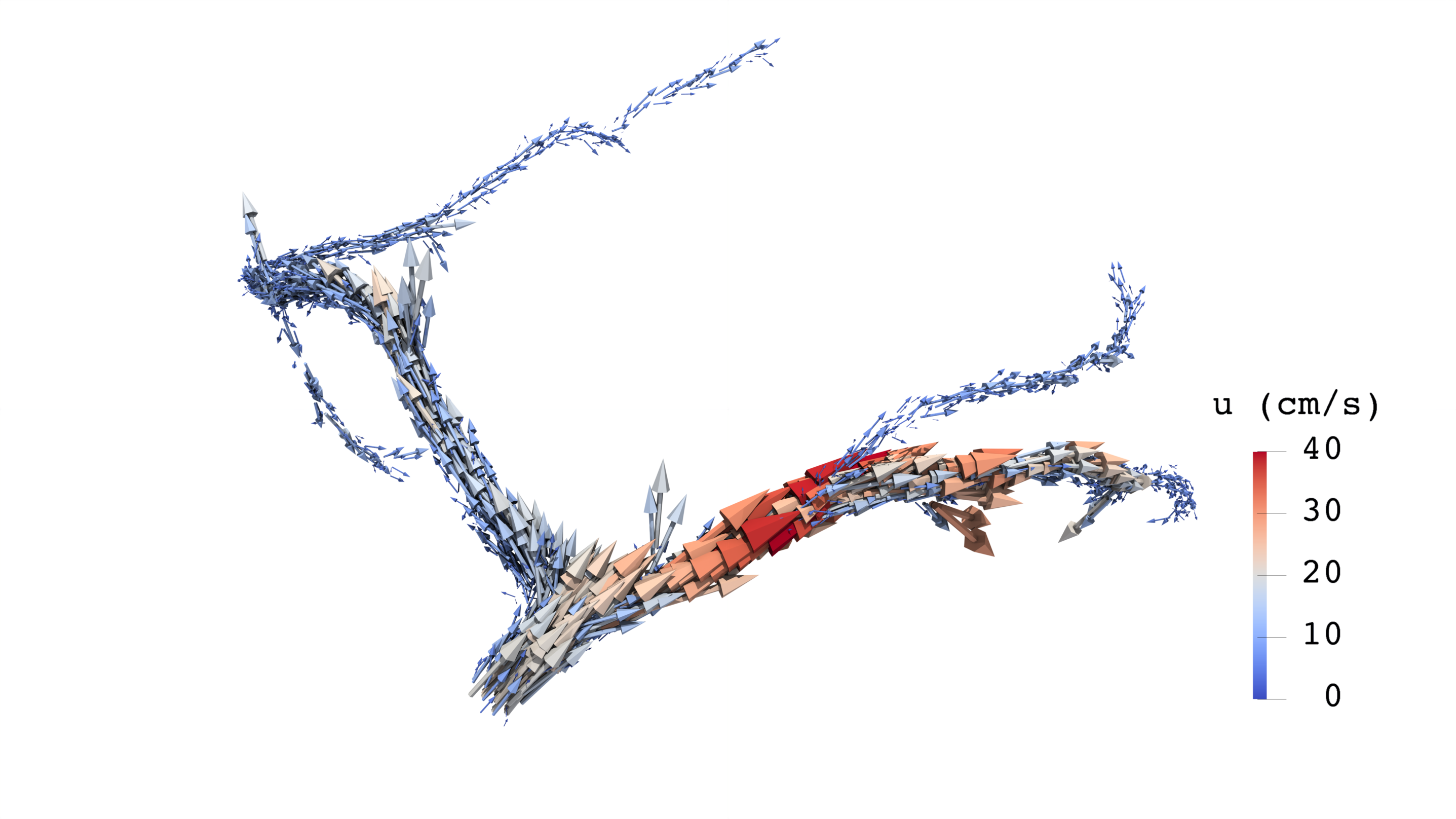} \label{fig:4d_flow_meas}}
    \caption{Measurement generation for the parameter estimation test cases.}
    \label{fig:4d_flow}
\end{figure}

To simulate 4D Flow MRI measurements, we followed a procedure similar to that
described in our previous work \cite{garay2022parameter}. Specifically, the same
CT velocity field used for the \emph{high-fidelity} dataset was first spatially
undersampled onto an image-like tetrahedral mesh with a resolution of $1
\ \text{mm}^3$ (see Figure \ref{fig:4d_flow_mesh}). The mesh was generated using
the algorithm reported in \cite{NOLTE2021102195}. Then, a complex magnetization
field was produced by perturbing the interpolated velocity with Gaussian noise at
$22 \ \text{dB}$. Finally, the velocity was reconstructed from the magnetization phase
using a velocity encoding parameter set to $120\%$ of the maximum velocity, in
order to avoid velocity aliasing. The final result is depicted at peak velocity
($t=0.69 \ \text{s}$) in Figure \ref{fig:4d_flow_meas}. Throughout the remainder of
this article, we refer to this as the \emph{Flow MRI-like} dataset.


\subsection{Inverse problem setup}

We define three test cases with an increasing complexity.

\begin{itemize}
    \item \textbf{Case 1.} We estimate four parameters,  $\ell_{5},
        \ell_{7},\ell_{9},\ell_{12}$,  corresponding to the outlets furthest from the
        inlet, as shown in Figure \ref{fig:coronary_mesh}. All other parameters are
        kept fixed during the estimation. 
        For the measurements, we consider first the \emph{high-fidelity} dataset,
        in order to obtain the best possible estimation we can get with this method.
    \item \textbf{Case 2.} The same four parameters are estimated from the
        \emph{4DFlow MRI-like} generated as described in Section \ref{sec:measurement_gen}.
    \item \textbf{Case 3.}  We estimate 16 out of 17 model parameters, namely
        $\ell_1,\dots,\ell_{16}$, from the \emph{high-fidelity} measurements.  The
        reason for fixing a single parameter is that the inverse problem requires a
        known pressure level; otherwise, the pressure field would only be determined
        up to an additive constant.
\end{itemize}

In all three cases, the initial guesses for the parameters were the mean value of the
reference duct lengths reported in Table \ref{tab:Parameters}, which was equals to $
\ell^0_m = 2.8 \rm{ mm}$. In order to ensure positivity of the estimated parameters, a
reparametrization was performed on the estimated lengths $\ell_m$ of the form $\ell_m
= \ell^0_m 2^{\beta_m}$, where the ROUKF method now optimizes for the $\beta_m$.
Moreover, the initial standard deviation for the estimation of $\beta_m$ was set to $0.5$.

\subsection{Results}

\begin{figure}[htbp]
    \centering
    \subfloat[Case 1: High-fidelity measurements]{
    \includegraphics[trim=0 0 0 10, clip, width=0.7\textwidth]{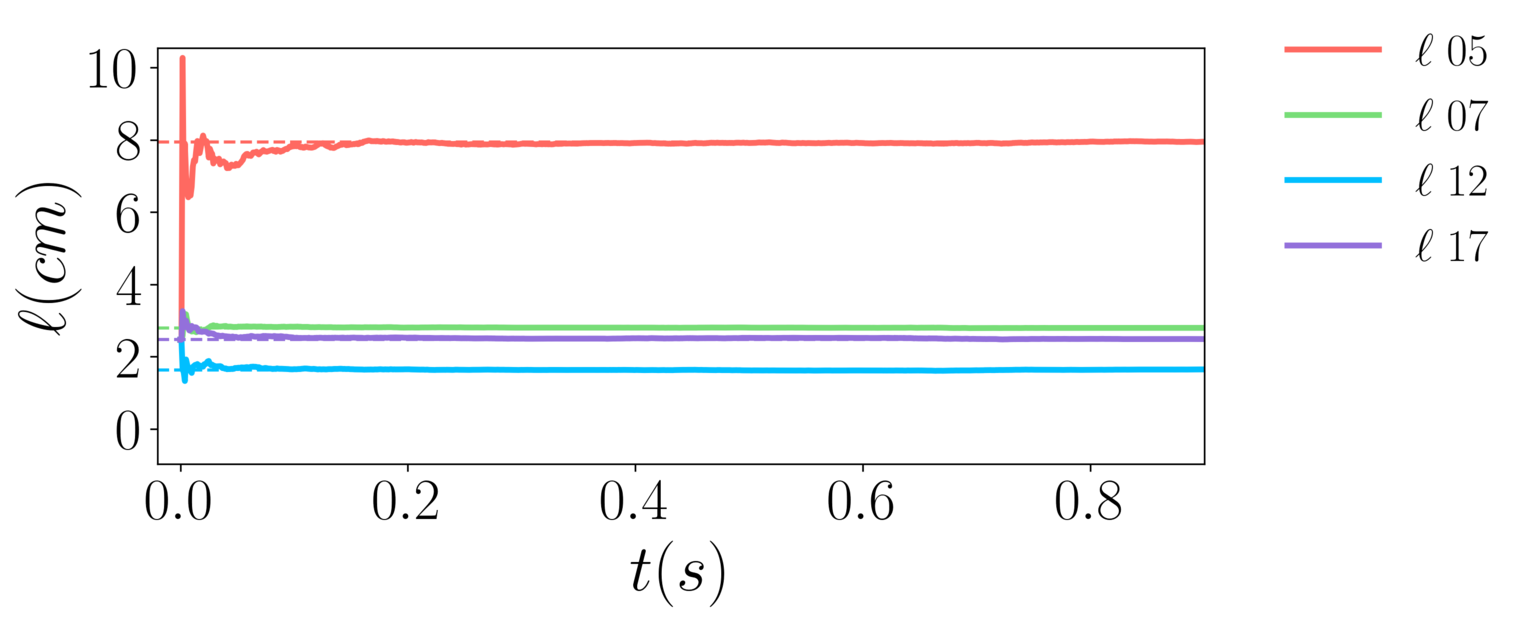} }\\
    \subfloat[Case 2: 4D Flow-like measurements]{
        \includegraphics[trim=0 0 0 10, clip,
    width=0.7\textwidth]{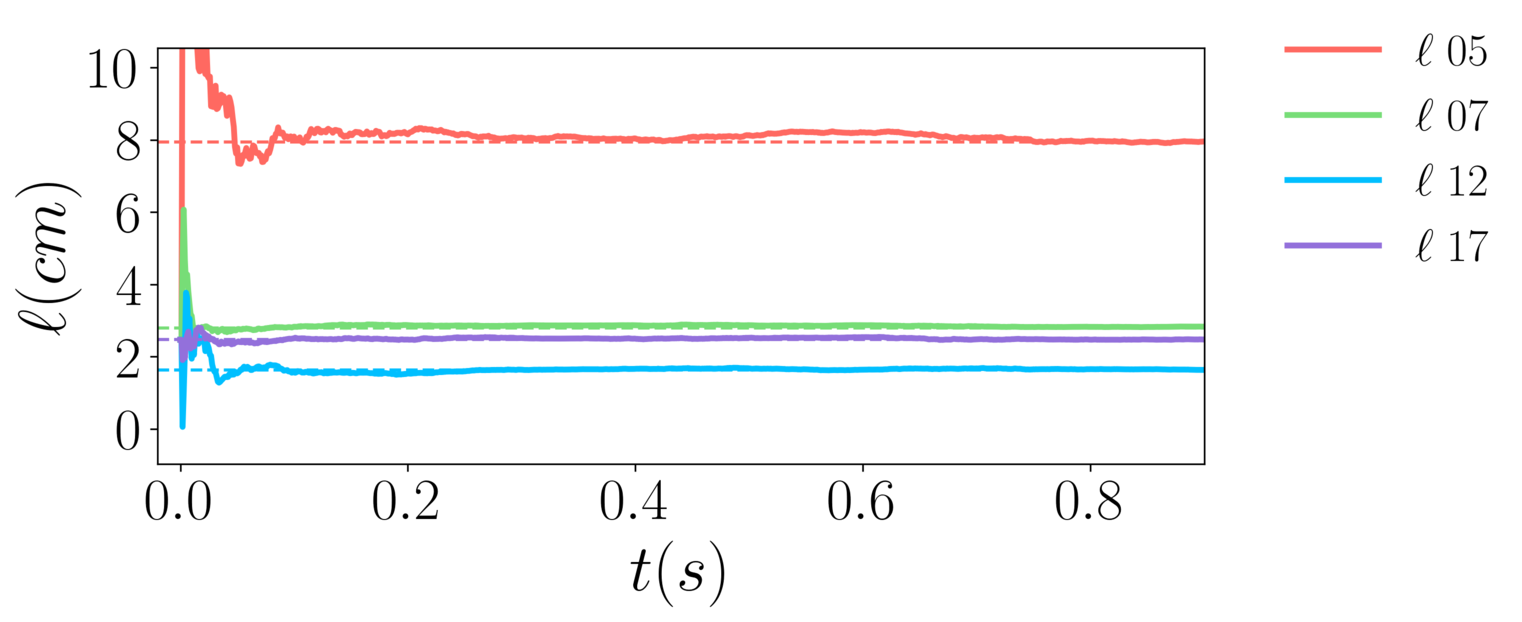} }
    \caption{Parameter estimation of 4 outlet lengths (Cases 1 and 2).}
    \label{fig:KalmanP1P2_p04}
\end{figure}

Figure \ref{fig:KalmanP1P2_p04} shows the estimation on only the four selected outlets
when using (Cases 1 and 2). Furthermore, Table \ref{tab:kalman_results_01} shows the
final estimated parameters in which it can be observed that when using the
high-fidelity data, the total mean error of the estimation was around $0.44 \%$. When
using the simulated 4D Flow measurements, the mean error was slightly raised to $0.48\%$.

\begin{table}[htbp] \centering
    \begin{tabular}{cc|cc|cc|cc}
        & & \multicolumn{2}{c|}{\makecell{Initial guess}} &
        \multicolumn{2}{c|}{\makecell{High-fidelity data}} &
        \multicolumn{2}{c}{\makecell{4D Flow-like data}} \\
        \midrule
        Boundary & $\ell_{ref}$ & $\ell_{init}$ & $\epsilon(\ell_{init})$ &  $\ell_{estim}$
        & $\epsilon(\ell_{estim})$  & $\ell_{estim}$ &
        $\epsilon(\ell_{estim})$ \\
        \midrule
        $\Gamma_{5}$       & {7.94} & {2.8} & {-64.74\%} & {7.92} &  {0.13\%} & {7.96}
        & {0.25\%}\\
        $\Gamma_{7}$       & {2.80} & {2.8} & {0.00\%} & {2.80} &  {0.00\%} & {2.83} &
        {1.07\%} \\
        $\Gamma_{12}$      & {1.63} & {2.8} & {71.78\%} & {1.65} &  {1.23\%} & {1.64}
        & {0.61\%} \\
        $\Gamma_{17}$      & {2.48} & {2.8} & {12.90\%} & {2.49} &  {0.40\%} & {2.48} & {0.00\%}
    \end{tabular}
    \caption{ROUKF estimations when using the high-fidelity and the 4D Flow-like
        measurement sets. The error was computed as:
        $\epsilon(\ell) = ( \ell -
        \ell_{ref})/\ell_{ref}$, and $\ell_{init}$, $\ell_{estim}$ and $\ell_{ref}$ are the
    initial guess, estimated and reference values for the duct lengths, respectively.}
    \label{tab:kalman_results_01}
\end{table}

\begin{figure}[!hbtp]
    \centering
    \subfloat[$\Gamma_{1} - \Gamma_{4}$]{
        \includegraphics[trim=0 5 0 15, clip,
    width=0.65\textwidth]{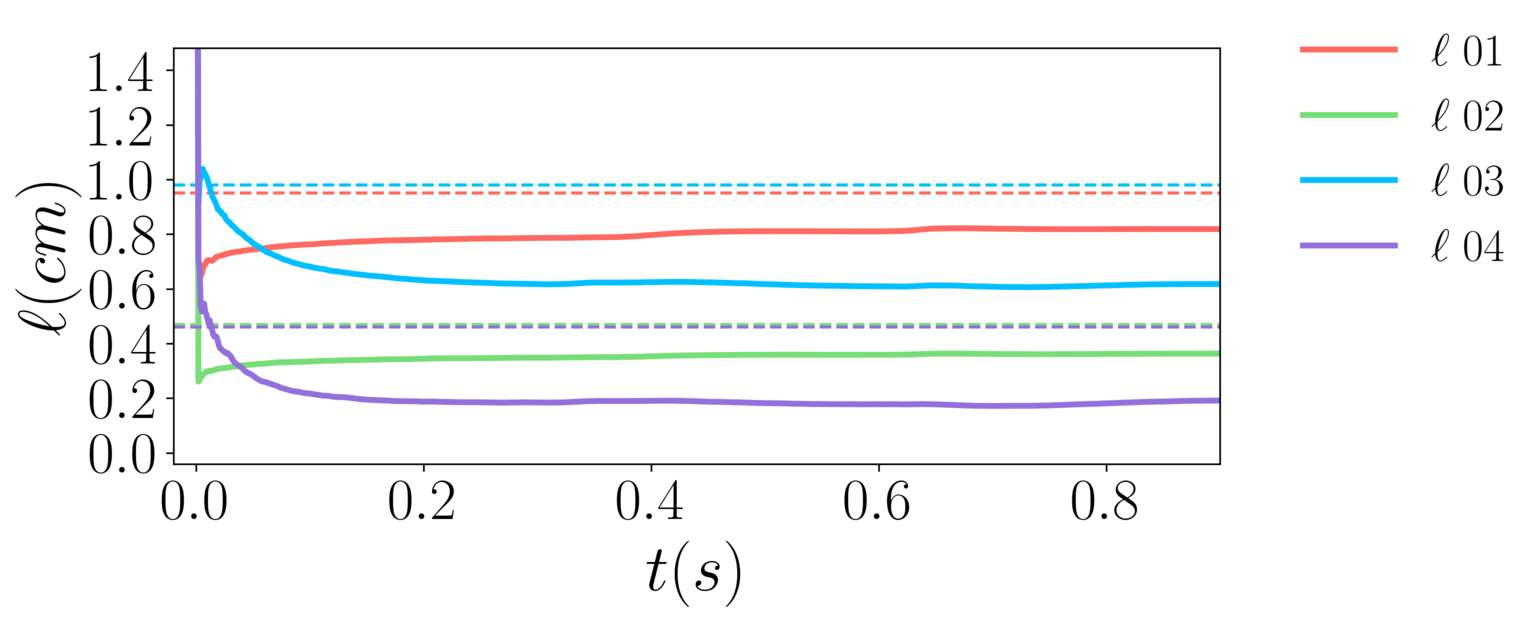} } \\
    \subfloat[$\Gamma_{5} - \Gamma_{8}$]{
        \includegraphics[trim=0 5 0 15, clip,
    width=0.65\textwidth]{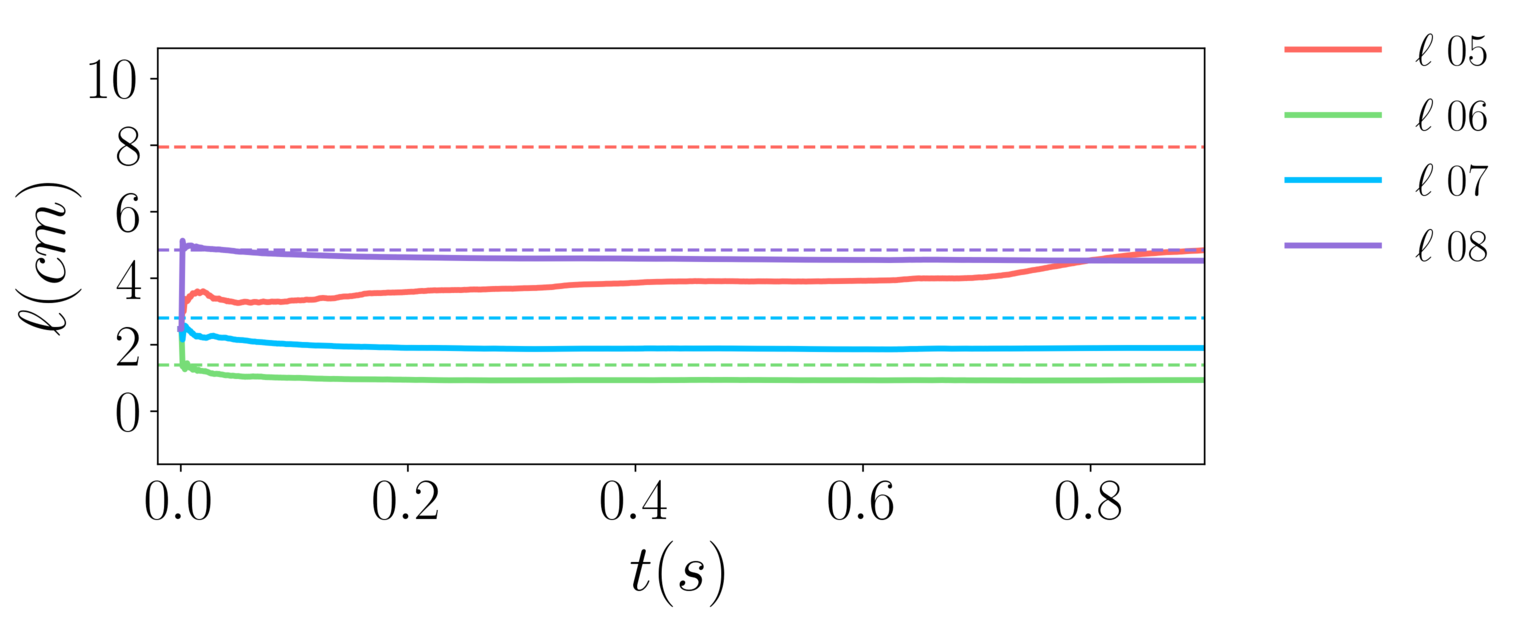} } \\
    \subfloat[$\Gamma_{9} - \Gamma_{12}$]{
        \includegraphics[trim=0 5 0 15, clip,
    width=0.65\textwidth]{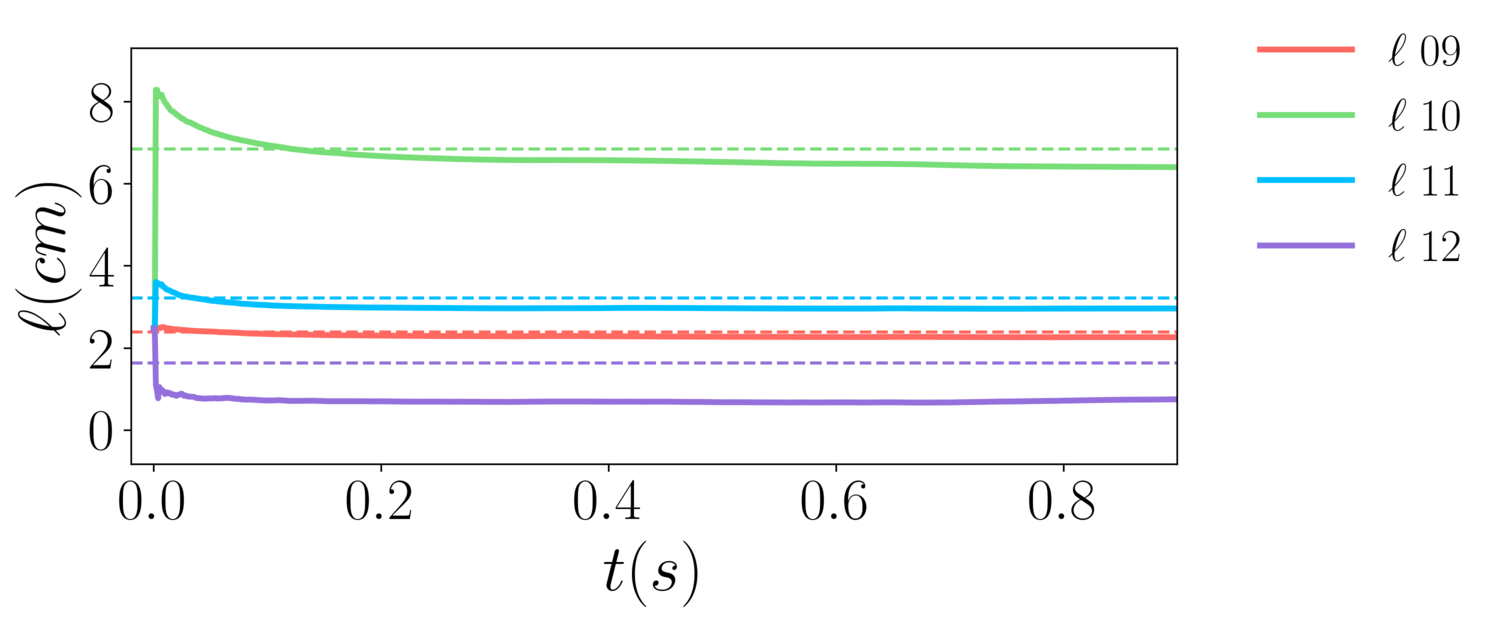} } \\
    \hspace*{-0.5cm}\subfloat[$\Gamma_{13} - \Gamma_{16}$]{
        \includegraphics[trim=0 5 0 15, clip,
    width=0.65\textwidth]{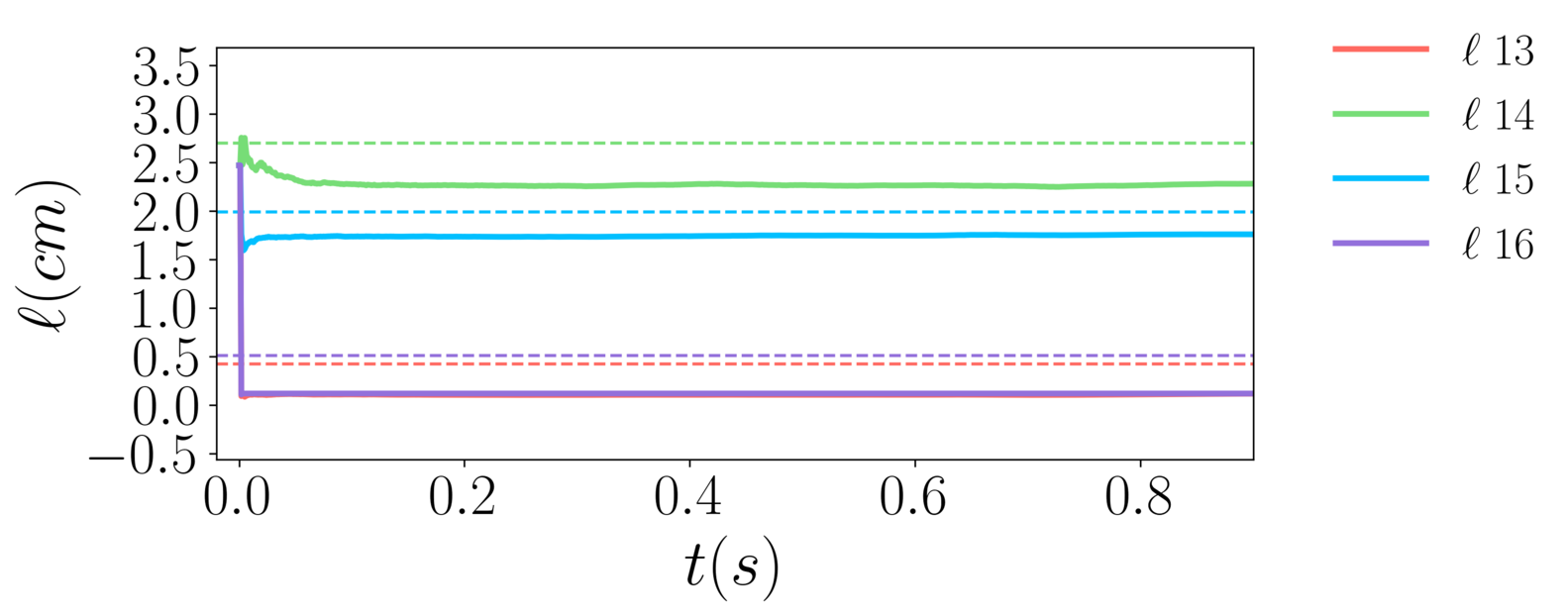} }
    \caption{Parameter evolution during the ROUKF run. Continuous lines shows the
    estimated parameter value while dashed lines are the reference values.}
    \label{fig:KalmanP1_all}
\end{figure}

Figure \ref{fig:KalmanP1_all} shows the ROUKF results for Case 3. It can be seen
that some parameters converge to stable values more rapidly than others, potentially
highlighting differences in identifiability.

Based on the initial guess, the average relative error across all parameters was
approximately $153\%$.  After the ROUKF run, this error was reduced to
around $30.77\%$. A summary of the final estimated values and their relative errors is
shown in Table \ref{tab:kalman_results_02}. These results demonstrate that the DuBC
method remains stable when varying the parameters as done by ROUKF, confirming its
suitability for parameter estimation problems.
It is worth noting that all estimated values remained systematically underestimated. A
possible explanation is as follows: since most of the initial guesses are smaller than
the target values for many of the parameters, the regularization imposed by the ROUKF
constrains their ability to reach the target values. Consequently, parameters whose
target values are below the initial guesses must also be reduced in order to produce
the target flow split.

\begin{figure}[!hbtp]
    \centering
    \subfloat[]{
        \includegraphics[trim=0 0 0 0, clip,
    width=0.48\textwidth]{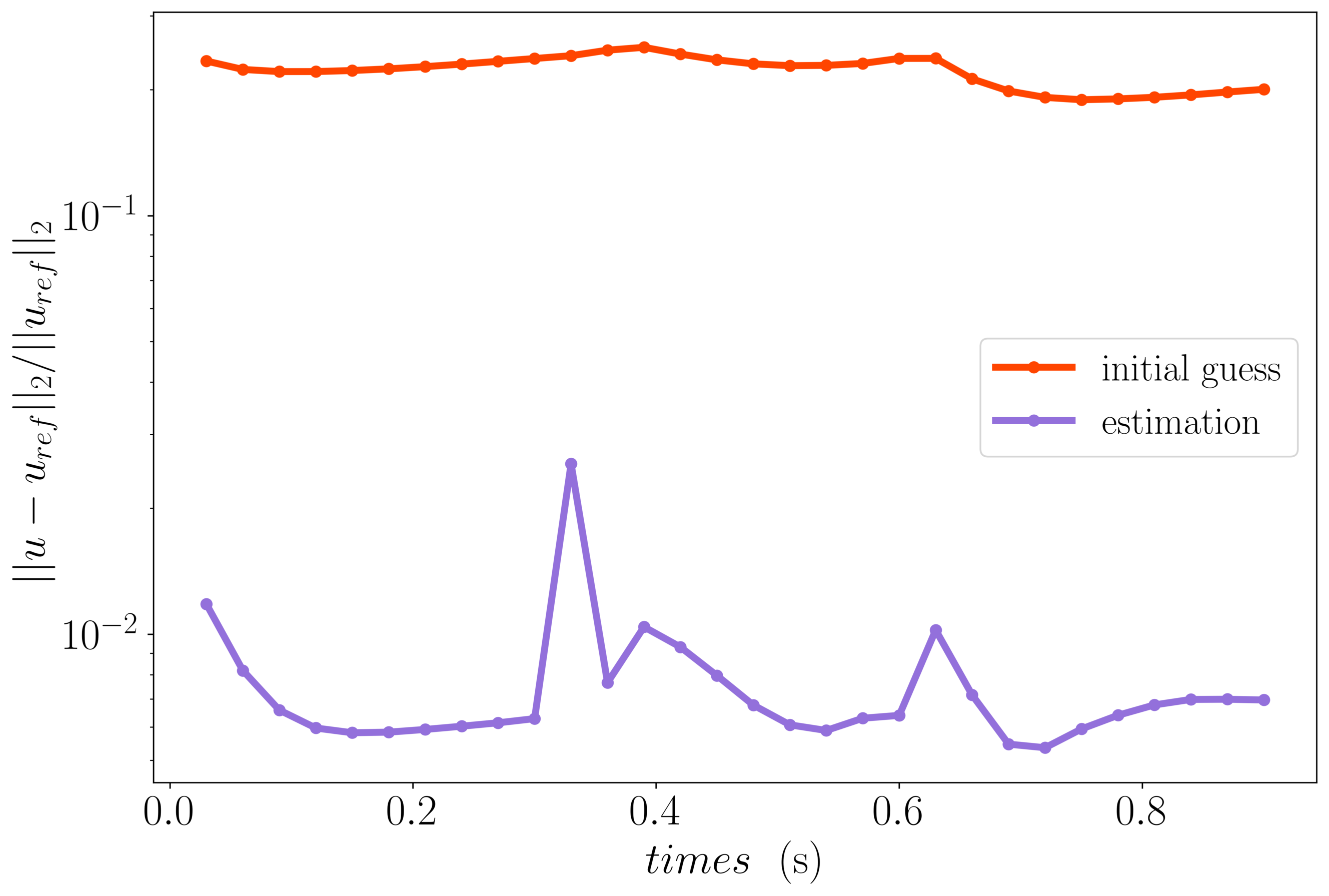} }
    \subfloat[]{
        \includegraphics[trim=0 0 0 0, clip,
        width=0.48\textwidth]{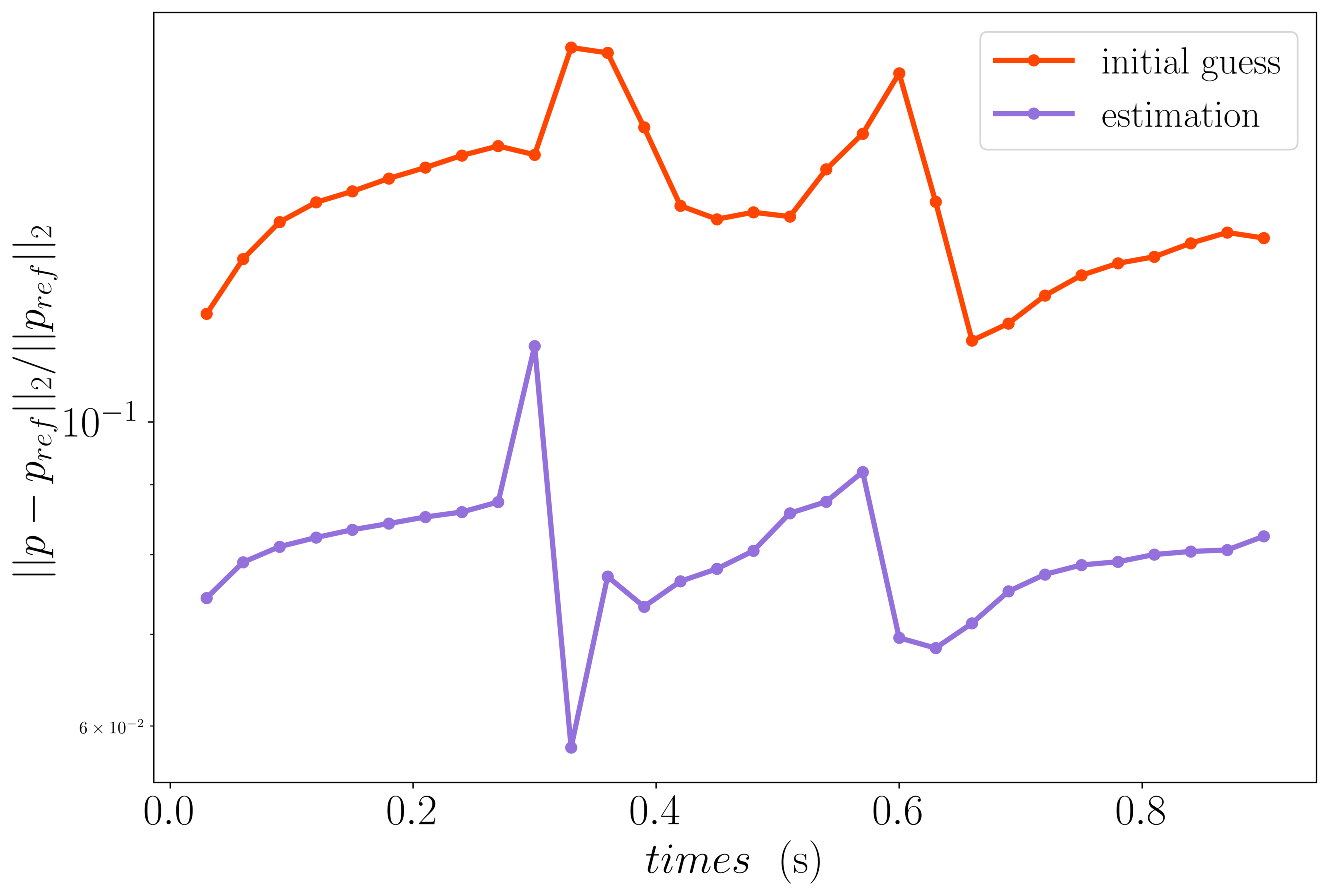}
    }
    \caption{Relative $L_2$ norms of the velocity (a) and pressure (b) fields obtained
        with the true set of parameters against the initial guesses and estimated
        parameters set for the Case 3. In both cases, the ROUKF method produced a
    reduction of the error over time.}
    \label{fig:kalman_results_norms}
\end{figure}

Finally, Figure \ref{fig:kalman_results_norms} shows the relative norm of the
velocity and pressure difference fields, between the true parameter solution, and
the solution obtained when we used the initial guess parameters as well as the
final estimated parameters for the Case 3 estimation case. From these curves, it
is evident that the parameter estimation significantly reduced the velocity error
by nearly an order of magnitude. The pressure error was also reduced, although to
a lesser extent.

\begin{table}[htbp] \centering
    \begin{tabular}{ccccc}
        Boundary & $\ell_{ref}$ & $\ell_{estim}$ & $\epsilon(\ell_{init})$ &
        $\epsilon(\ell_{estim})$ \\
        \midrule
        $\Gamma_{1}$ & {0.95} &  {0.82} & {194.7\%} & {-13.68\%} \\
        $\Gamma_{2}$ & {0.47} &  {0.36} &  {495.7\%} & {-23.40\%} \\
        $\Gamma_{3}$ & {0.98} &  {0.62} & {185.7\%} & {-36.73\%} \\
        $\Gamma_{4}$ & {0.46} &  {0.19} & {508.7\%} & {-58.70\%} \\
        $\Gamma_{5}$ & {7.94} &  {4.84} & {-64.7\%} & {-39.04\%} \\
        $\Gamma_{6}$ & {1.38} &  {0.93} & {102.9\%} & {-32.61\%} \\
        $\Gamma_{7}$ & {2.80} &  {1.90} & {0.0\%} & {-32.14\%} \\
        $\Gamma_{8}$ & {4.85} &  {4.52} & {-42.27\%} & {-6.80\%} \\
        $\Gamma_{9}$ & {2.39} &  {2.26} & {17.15\%} & {-5.44\%} \\
        $\Gamma_{10}$ & {6.84} & {6.40} & {-59.06\%} & {-6.43\%} \\
        $\Gamma_{11}$ & {3.21} &  {2.96} & {-12.77\%} & {-7.79\%} \\
        $\Gamma_{12}$ & {1.63} &  {0.74} & {71.78\%} & {-54.60\%} \\
        $\Gamma_{13}$ & {0.42} &  {0.12} & {566.67\%} & {-71.43\%} \\
        $\Gamma_{14}$ & {2.70} &  {2.28} & {3.70\%} & {-15.56\%} \\
        $\Gamma_{15}$ & {1.99} &  {1.76} & {40.70\%} & {-11.56\%} \\
        $\Gamma_{16}$ & {0.51} &  {0.12} & {449.0\%} & {-76.47\%}
    \end{tabular}
    \caption{ROUKF estimations when using the high-fidelity
        measurements for 16 out of 17 parameters. The error was computed as:
        $\epsilon(\ell) = ( \ell -
        \ell_{ref})/\ell_{ref}$, and $\ell_{init} = 2.8$, $\ell_{estim}$ and $\ell_{ref}$ are the
    initial guess, estimated and reference values for the duct lengths, respectively.}
    \label{tab:kalman_results_02}
\end{table}

\section{Conclusions}
\label{sec:conclusion}

In this work, we presented a new duct boundary condition (DuBC) and
demonstrated its application to the simulation of coronary flows. This boundary
condition is a special implementation of the method of asymptotic partial decomposition
of a domain (MAPDD). We provided its extension to a fractional step scheme and tested it
on an inverse problem in a complex hemodynamic setting.

The simplicity of implementation and stability properties of DuBC make it appealing for
applications involving highly ramified domains such as coronary or cerebral
arteries. However, this comes at the cost of introducing "virtual" distal vasculature
lengths at the domain boundaries—a parameter that is not directly measurable and
may introduce uncertainty if not chosen carefully.

While this approach simplifies the problem compared to Windkessel-type boundary
conditions, which require estimating several lumped parameters, it does not benefit
from the availability of well-established physiological reference values in the
literature. Nevertheless, we demonstrate that the virtual lengths can be estimated
from velocity data, in a manner analogous to how resistances are derived in Windkessel models.

A natural next step for this framework is its application to real 4D Flow MRI
data, as reported in \cite{blanken2021coronary}. Once the DuBC parameters are
calibrated, forward simulation with DuBC can be used to reconstruct hemodynamic
fields that are consistent with the measurements, enabling subject-specific blood
flow analysis in branching vessel domains.



\section*{Acknowledgements}
\label{sec:ack}

C.B. and D.N. acknowledge the funding from the European Research Council (ERC) under
the European Union's Horizon 2020 research and innovation programme (grant agreement
No 852544 - CardioZoom). J.G. acknowledges the funding of ANID Chile by the Fondecyt
Postdoc project No 3230549.


\bibliography{biblio_mapdd.bib}

\end{document}